\documentclass[11pt]{article}

\usepackage{bbding}
\usepackage{pifont}
\usepackage{amsfonts}
\usepackage{graphicx}
\usepackage{amsmath}
\usepackage{amsthm}
\usepackage{amssymb}
\usepackage{concrete}

\def\beq{\begin{equation}}
\def\eeq{\end{equation}}
\def\nd{\noindent}

\def\<{\leq}
\def\>{\geq}

\newtheorem{thm}{Theorem}[section]
\newtheorem{lem}{Lemma}[section]
\newtheorem{prop}{Proposition}[section]

\newtheorem{defi}{Definition}[section]
\newtheorem{rem}{Remark}[section]

\begin{document}
\title{Algebras associated with Pseudo Reflection Groups: A Generalization of Brauer Algebras }
\author{Zhi Chen}
%\address{Department of Mathematics, University of science and Technology of China, Hefei, Anhui province, 230026, China} \email{zzzchen@ustc.edu.cn}
\maketitle

\begin{abstract}
\noindent We present a way to associate an algebra $B_G (\Upsilon) $
with every pseudo reflection group $G$. When $G$ is a Coxeter group
of simply-laced type we show $B_G (\Upsilon)$ is isomorphic to the
generalized Brauer algebra of simply-laced type introduced by
Cohen,Gijsbers and Wales[10]. We prove $B_G (\Upsilon )$ has a
cellular structure and be semisimple for generic parameters when $G$
is a rank 2 Coxeter group. In the process of construction we
introduce a Cherednik type connection for BMW algebras and a
generalization of Lawrence-Krammer representation to complex braid
groups associated with all pseudo reflection groups.
\end{abstract}

%%%%%%%%%%%%%%%%%%%%%%%%  PRELIMINARIES %%%%%%%%%%%%%%%%%%%%%%%%%%%%%
\section{Introduction}

The Brauer algebras $B_n (\tau )$ were introduced by Brauer in
1930's, motivated by the purpose of creating a Schur-Weyl duality
for orthogonal groups and symplectic groups. Roughly speaking the
Brauer algebra $B_n (\tau )$ is an algebra containing the group
algebra of symmetric group $k S_n $, and sharing many important
 properties with $kS_n $.  First, $B_n (\tau )$ is semisimple for
generic $\tau$ as proved in [27], whence their irreducible
representations are also labeled by Young diagrams.  Secondly just
as $k S_n $ has Hecke algebra as a deformation , which provide a way
to define the HOMFLY polynomial invariants for Links, Brauer
algebras have a natural deformation named Birman-Murakami-Wenzl
algebras(also called BMW algebras),which dominate the Kauffman
polynomial invariants in similar way [3][23]. The Brauer algebras
also admit cellular structures like $kS_n $ in the sense of Graham
and Lehrer[15].

As many objects related to Lie theory, In [10] Cohen, Gijsbers and
Wales defined BMW algebras and Brauer algebras for simply-laced
Coxeter groups. Those new algebras have almost all important
algebraic properties of $B_n (\tau )$. They were shown to be
semisimple for $D_n $ type for generic parameters, and admitting
cellular structures ([11]). In [10] the author ask whether there
exist a definition of generalized BMW algebras for non-simply laced
Coxeter groups.  In another way in [16] Haring-Oldenburg  defined
Cyclotomic BMW algebras and their degenerate version Cyclotomic
Brauer algebras. These algebras were studied in many subsequent
papers.

In this paper we introduce an algebra $B_G (\Upsilon )$ for every
Coxeter group and every pseudo reflection group $G$, where
$\Upsilon$ is a set of parameters.
 These algebras have the following properties .

 If $G$ is finite then $B_G (\Upsilon )$
is a finite dimensional algebra containing the group ring $kG $
naturally.

There exist a flat formal connection which can deform every
representation $B_G (\Upsilon )$ to a one-parameter family of
representations of $A_G$, the braid group associated with $G$. One
of these representations is the generalized Lawrence-Krammer
representation.

When $G$ is a Coxeter group of simply-laced type, $B_G(\Upsilon)$ is
isomorphic to the generalized Brauer algebra introduced in [10].
When $G$ is a pseudo reflection group of $G(m,1,n)$ type, the
cyclotomic Brauer algebra introduced in [16] is isomorphic to a
natural quotient of our algebra $B_G(\Upsilon )$.

 When $G$ is a rank 2 Coxeter group the algebra $B_G(\Upsilon )$ is
semisimple for generic $\Upsilon $ and admit a cellular
structure(Thm 5.2,Prop5.2 ).

Two major ingredients in the construction of $B_G(\Upsilon)$ are
Cherednik type connections for BMW algebras and the generalized
Lawrence-Krammer representations of Cohen,Wales [9] and Marin [22].

Let $H_n (q)$ be the Hecke algebra of the symmetric group $S_n $.
Let
$$Y_n =\{(z_1 ,\cdots ,z_n )\in \mathbb{C} ^n | z_i
\neq z_j ,\ for\ any\ i\neq j \}.$$ It is well-known that for
generic $q$ the set of irreducible representations of $H_n (q)$ is
in one to one correspondence with the set of irreducible
representations of $S_n $. This correspondence has a geometric
description. The Cherednik connection describe how a representation
of $S_n $ can be deformed to a representation of $H_n (q)$.

 Let $\omega _{i,j}= d\log (z_i -z_j )$ for $i<j$ be closed forms on $Y_n $.
The Cherednik connection is a formal flat connection
  $$\Omega _n = \kappa \sum _{i<j} s_{i,j} \omega _{i,j}$$

 defined on the bundle $Y_n \times \mathbb{C} S_n $.  Where $s_{i,j}\in S_n $ is the $(i,j)$ permutation.

Cyclotomic Hecke algebras also have such flat connections as shown
in [5].  Let $G$ be a pseudo reflection group acting on a complex
linear space $E$, let $H_G (\bar{q})$ be the associated cyclotomic
Hecke algebra. We denote the set of reflection hyperplanes of $G$ as
$\{ H_v \} _{v\in P} $, and denote the set of pseudo reflections of
$G$ as $R$. Define $v: R\rightarrow P $ by requesting the reflection
hyperplane of $s$ to be $H_{v(s)}$. Denote $M_G = E\setminus \cup
_{i\in P} H_i $ being the complementary space of reflection
hyperplanes. By Steinberg theorem [26] the induced action of $G$ on
$M_G $ is free. Complex braid group $A_G$ associated with $G$ is
defined as $\pi _1 (M_G /G)$. For any $v\in P$ let $f_v $ be any
nonzero linear function on $E$ with kernel $H_v $, let $\omega _v =
d\log f_v $. The Cherednik type connection of $G$ is defined on the
bundle $M_G \times \mathbb{C} G $ as follows

$$\Omega _G  = \sum _{v\in P} \kappa ( \sum _{s\in R, v(s)=v }
\mu _s s ) \omega _v .$$

Where $\{ \mu _s \} _{s\in R} $ is a set of parameters satisfying
$\mu _{s_1 } = \mu _{s_2 } $ if $s_1 $ is conjugate to $s_2$.
Suppose $(V,\rho )$ is linear representation of $G$, then we have a
flat connection on bundle $M_G \times V$:

$$\rho ( \Omega _G  ) = \sum _{v\in P} \kappa ( \sum _{s\in R,
v(s)=v } \mu _s \rho ( s) ) \omega _v .$$

$\rho ( \Omega _G  )$ induce a flat connection on the quotient
bundle $M_G \times _G V$ which define a one parameter class of
representations of $A_G$ through monodromy. These monodromy
representations factor through $H_G (\bar{q})$ for suitable
$\bar{q}$.[5]

 It is known that for generic parameters the
set of irreducible representations of BMW algebras is also in one to
one correspondence with the set of irreducible representations of
the Brauer algebra $B_n (\tau )$.  It is reasonable to expect that
there exist Cherednik type connection for BMW algebras.

We find the formal connection $\bar{\Omega } _{n} = \kappa \sum
_{i<j}( s_{i,j} -e_{i,j} ) \omega _{i,j}$ defined on $Y_n \times B_n
(\tau )$ is just what we want. It is flat and can deform every
representation of Brauer algebra to a representation of BMW
algebra(Thm 3.2 ).  Where $e_{i,j} \in B_n (\tau )$ is the element
described by the following graph.

\begin{figure}[htbp]

  \centering
  \includegraphics[height=4.5cm]{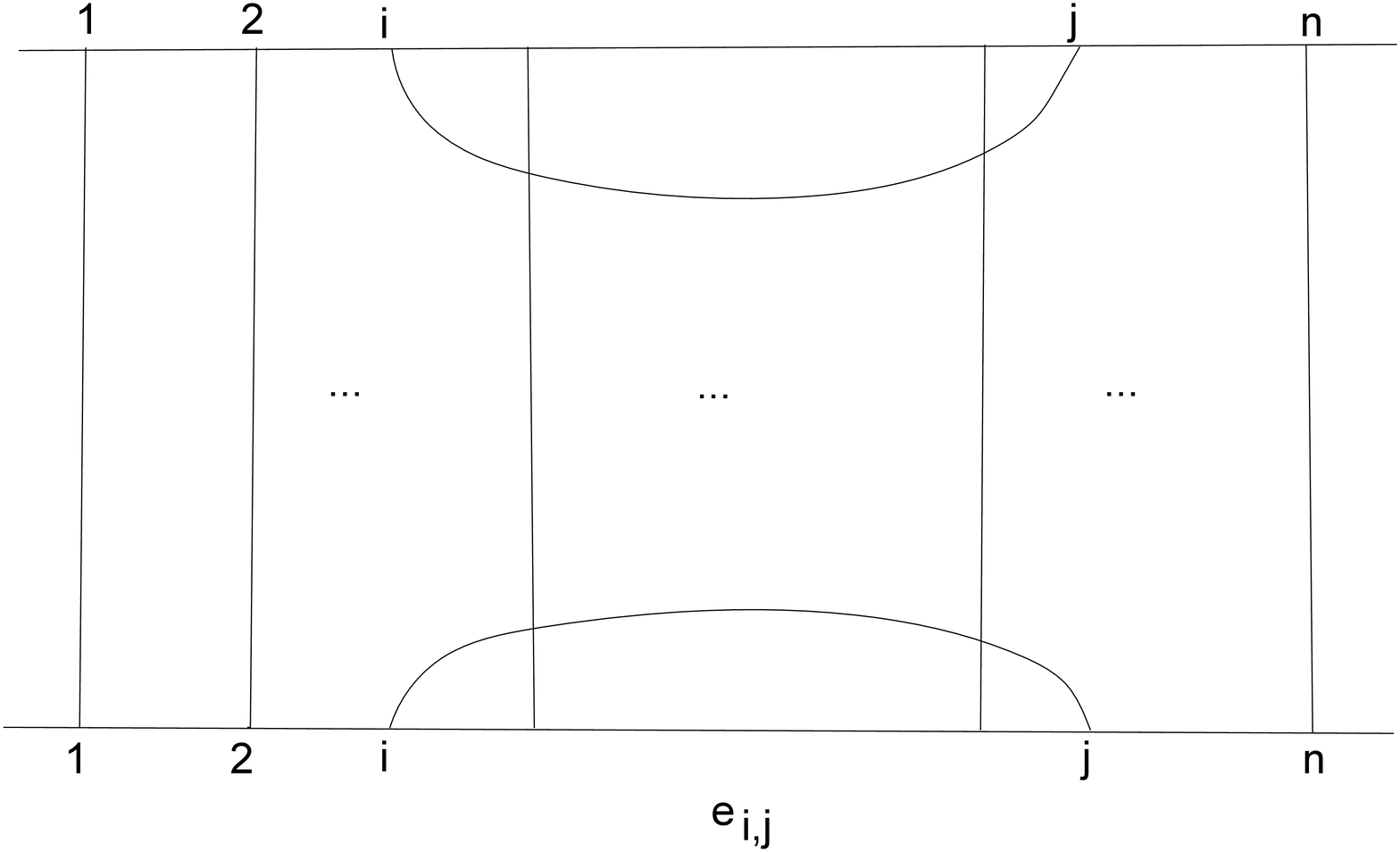}

\end{figure}

Now we see the deformable property  of  $\mathbb{C} G$ and $B_n
(\tau )$  embody in the existence of a flat formal connection
describing the deformation. We hope generalized Brauer algebras
should be deformable and have such flat connections as well. The
concise form of above flat connection motivate us to imagine that
$B_G (\Upsilon )$ should be generated by elements in $G$ along with
a set of special
 elements $\{ e_i \} _{i\in P}$, such that the following formal
 connection $\bar{\Omega } _{G}$ on $M_G \times B_G (\Upsilon )$ is flat and $G$-equivariant.
$$\bar{\Omega } _{G} = \sum _{v\in P} \kappa (  \sum _{s\in R, v(s)=v } \mu _s
s  - e_v ) \omega _v .$$

If all pseudo reflections in $G$ have rank 2, and if we set all
parameters $\mu _s =1$, then the connection become: $\bar{\Omega }
_{G} = \sum _{v\in P} \kappa ( s_v   - e_v ) \omega _v $, which has
a form similar with $\bar{ \Omega } _n $.

The flatness condition give some relations between $s_i $'s and $e_i
$'s.  But these relations are too coarse to produce an interesting
finite dimensional algebra.

Annother common  feature of  BMW algebras of simply-laced types is
that they all contain the generalized Lawrence-Krammer
representations(from now on we call them by LK representations)
defined in [9](see[10]). Infinitsiml version of LK representations
are found by Marin[21]. They are described by some concise  flat
connections. Marin then generalized these flat connections to all
complex reflection groups( pseudo reflection groups whose pseudo
reflections all have rank 2), thus generalize LK representations to
all complex reflection groups [22](see 5.2). In section 4 We show LK
representations can actually be generalized to all pseudo reflection
groups and have a very simple characterization.

Let $V_G = \mathbb{C} \{ v_i \} _{i\in P }$ be a vector space with a
basis in 1 to 1 correspondence with $\{H_i \} _{i\in P}$, the set of
reflection hyperplanes of $G$. Since $G$ permutes the set of
reflection hyperplanes, there is a natural
representation $\iota : G \rightarrow GL(V_G )$.\\
\nd { \bf Theorem 4.2 } {\it   For every $i\in P$, suppose $p_i \in
End(V_G )$ is a projection to $\mathbb{C} v_i$. Explicitly suppose
$p_i (v_j ) = \alpha _{i,j} v_i $ for $j\neq i$; $p_i (v_i ) =m_i
v_i $.  Then the following connection on $M_G \times V_G $ , $\Omega
_{LK} = \sum _{i\in P } \kappa (\sum _{s:v(s) =i } \mu _s \iota (s)
- p_i ) \omega _i$ is flat and $G$-equivariant if and only if
$\alpha _{i,j} = \sum _{s: \iota(s)(v_j ) =v_i} \mu _s $,  and $m_i
= m_j  $ if there exist $w\in G$ such that $w(H_i ) = H_j $.  }
\\

When $\Omega _{LK}$ satisfy the $G$ invariant and flatness condition
in Theorem 4.2, it induce a flat connection $\bar{\Omega}_{LK}$ on
the quotient bundle $M_G \times _G V_G$. We define the generalized
LK representation of $A_G$ as the monodromy representation of
$\bar{\Omega}_{LK}$. When $\mu _s =1$ for all $s$, the connection
$\Omega _{LK} $ become the flat connection of Marin [22].

To find more relations for $B_G (\Upsilon )$  now we make annother
request that all generalized BMW algebras should contain generalized
LK representations.  This condition can be presented as: there exist
some representation $(V,\ \rho )$ of $B_G (\Upsilon )$, such that
the connection $\rho ( \bar{\Omega} _{G } )$ coincide with $\Omega
_{LK}$. This request give us both restriction and indication on
relations we are looking for. We find the algebra defined in
Definition 5.1 satisfies
 these requests perfectly.

\nd {\bf Definition 5.1 } {\it  The algebra $B_G (\Upsilon )$
associated with pseudo reflection group $(V,G)$ is generated by
  $\{ \bar {w} \}_{w\in G } \cup  \{ e_i \}_{i\in P}$, submitting to the following
relations. Here $P$ is the index set of reflection hyperplanes of
$G$ as defined before.

\begin{itemize}
 \item[(0)]  $\bar {w_1} \bar {w_2} = \bar {w_3 }$ if $w_1 w_2 =w_3 $;
\item[(1)] $ \bar{s}_i  e_i  = e_i \bar{s}_i = e_i $, for $i\in P$;
\item[(2)] $e_i ^2 = m_i e_i $ ;
\item[(3)] $\bar {w} e_j  = e_i \bar {w} $ , if $w \in G  $ satisfies $w(H_j )=H_i$;
\item[(4)] $e_i e_j = e_j e_i $,
 if  $\{ \ k\in P |\  H_k \supset H_i \cap H_j \}
= \{i,j \}$ ;
 \item[(5)] $ e_i e_j = (\sum _{s\in R(i,j) } \mu _s s )e_j $ ,
 if $\{ k\in P | H_k \supset H_i \cap H_j \}   \neq
  \{i,j \},$
  and $R(i,j) \neq \emptyset $.
\item[(6)] $e_i e_j =0 $, if $\{  k\in I | H_k \supset H_i \cap H_j  \}
\neq
  \{i,j \} $ and $R(i,j) =\emptyset $.
\end{itemize}
where $R(i,j) = \{s\in R | s(H_j) =H_i \}$. $s_i $ is any pseudo
reflection fixing $H_i$.Constants $\{ \mu _s \} _{s\in R }$ satisfy
$\mu _{s_1 } =\mu _{s_2 } $ if $s_1 $ is conjugate to $s_2 $ .
Constants $\{ m_i \}_{1\in P } $ satisfy $m_i =m_j$ if there exists
 $w\in G $ such that $w(H_i ) =H_j $. The symbol $\Upsilon $ means data $\{ \mu _s \} \cup \{ m_i \}$.} \\

For these algebras $B_G (\Upsilon)$, we prove they are finite
dimensional( Thm 5.1 ), the connection $\bar{\Omega} _G $ is
$G$-invariant and flat (Prop 5.1). And It is readily to see it could
realize the connection $\Omega _{LK}$.  We hope these algebras
satisfy other important properties of Brauer algebra: Semisimple for
generic parameters;  have cellular structures; can be deformed to
certain generalized BMW algebras. Existence of the connection
$\Omega _{Br}$ seems supporting the last property.  Section 5.2 is
devoted to study of cases when $G$ are dihedral groups. We prove
they are semisimple for generic parameters(Thm 5.2,Thm 5.5), and
have cellular structures(Thm 5.3, Prop 5.2).  We hope these two
properties hold for all $B_G (\gamma )$.

In section 5.3 we  construct canonical presentations of $B_G
(\Upsilon )$ when $G$ is a Coxeter group of finite type (Thm 5.10),
when $G$ is of simply-laced type, it can be seen from this
presentation  that $B_G (\Upsilon )$ coincide with the generalized
Brauer algebra in [9].  When $G$ is a pseudo reflection group of
type $G(m,1,n)$ we also find a canonical presentation for $B_G
(\Upsilon )$ (Thm 5.11 ), from which  we see the cyclotomic Brauer
algebra introduced in [12] is isomorphic to a quotient of $B_G
(\Upsilon )$ by letting some element $e_i $ to be zero.

In the last section 5.4 we discuss the uniqueness of $B_G (\Upsilon
)$. If we replace (6) in definition 5.1 with a milder relation
$$(6) ^{'} : e_i e_j =e_j e_i , if \{  k\in I | H_k \supset H_i \cap H_j  \}  \neq
  \{i,j \}  and R(i,j) =\emptyset  $$

  and keep other relations, the resulted algebra $\hat{B} _G (\Upsilon
  )$ satisfies all above mentioned properties of $B_G (\Upsilon )$
  except semisimplity and existence of a cellular structure.
  Proposition 5.4 shows when $G$ is a rank 2 Coxeter groups,  for generic $\Upsilon$ the relation
  $(6)^{'}$ degenerate to relation $(6)$.  As we know BMW algebra
  have the same dimension for all parameters, proposition 5.4 seems
  proving $\hat{B} _G (\Upsilon)$ isn't a good choice.
  \\

\nd {\bf Acknowledgements} I would like to thank Professor Toshitake
 Kohno for teaching me KZ equations, and thank Professor Sen Hu and Professor Bin Xu
 for many beneficial communications  . Also thank Professor Jie Wu for his
 invitation to NUS in Dec 2008, this work was partially done during
 that stay.  And thank Professor Hebing Rui for many valuable advices.

\section{Preliminaries }
  \subsection{Brauer algebras and BMW algebras }
  Brauer algebra $B_n (\tau )$ is a graphic algebra
 in the sense that it has a basis
consisting of elements presented by graphs, and the relations
between them can be described through graphs.  $B_n (\tau ) $ has a
canonical presentation with generators  $s_1 , \cdots ,s_{n-1 }$,
$e_1 ,\cdots ,e_{n-1 }$ and relations listed in table 1.  $\mathcal
{B}_n (\tau )$ has a natural deformation discovered by Birman,
Murakami, Wenzl which are now called BMW algebras [3][23]. These
algebras support a Markov trace which gives the Kauffman polynomial
invariants of Links.  We denote these BMW algebras as $B _n (\tau ,l
)$.  Where $l$ is a parameter of deformation.  There is $B_n (\tau)
\cong B_n (\tau , 1 ).$ We list generators and relations of $B_n
(\tau ) $ and $B _n (\tau ,l )$ in table 1 according to [9]. Where
$m=\frac {l-l^{-1}}{1-\tau}$.
\\

\begin{tabular}{|l|l|l|}

  % after \\: \hline or \cline{col1-col2} \cline{col3-col4} ...
\hline

   & $B_n (\tau )$  & $B _n (\tau,l)$ \\
\hline

  Generators    & $s_1 $,$\cdots $,$s_{n-1}$;$e_1$ ,$\cdots$ ,$e_{n-1}$ & $X_1$ , $\cdots $,$X_{n-1}$; $E_1$ ,$\cdots$ ,$E_{n-1}$ \\
\hline

  Relations   &  $s_i s_{i+1} s_i =s_{i+1} s_i s_{i+1}$    & $X_i X_{i+1} X_{i} =X_{i+1} X_{i} X_{i+1}$\\
&$for 1\leq i\leq n-2 $; &$for 1\leq i\leq n-2$  ;\\

& $s_i s_{i-1} e_i = e_{i-1} s_i s_{i-1}  $
& $X_i X_{i-1} E_i = E_{i-1} X_i X_{i-1} $\\
& $for 2\leq i\leq n-1$;  &$for 2\leq i\leq n-1$;\\

 & $s_i s_{i+1} e_i = e_{i+1} s_i s_{i+1}  $
& $X_i X_{i+1} E_i = E_{i+1} X_i X_{i+1} $\\
&$for 1\leq i\leq n-2$; &$for 1\leq i\leq n-2$;\\
& $s_i s_j =s_j s_i$ ,$|i-j|\geq 2 $;& $X_i X_j =X_j X_i ,|i-j| \geq 2$ ;\\
 &$s_i ^2 =1 ,for\  all\  i $;& $l(X_i ^2 +mX_i -1) =mE_i , for\  all\
i $;\\
 & $s_i e_i =e_i ,for\ all\ i $;&$ X_i E_i = l^{-1} E_i ,for\ all\ i $;\\
& $e_i s_{i+1} e_{i} = e_{i}, 1\leq i \leq n-2 $; &$ E_i X_{i+1}
E_{i} = l E_i ,1\leq i \leq n-2 $;
\\
&$ e_i s_{i-1} e_i =e_i , 2\leq i \leq n-1$ ; &$ E_i X_{i-1} E_{i} =
lE_{i} ,2\leq i\leq n-1 $; \\
&$s_i e_j = e_j s_i , |i-j|>1 $; & $X_i E_j =E_j X_i ,|i-j|>1 $;\\
&$e_i ^2 = \tau e_i ,for\  all\  i.$ &$ E_i ^2 =\tau E_i $
. \\

  \hline
\end{tabular}
\\
\\

The structure of Brauer algebras and BMW algebras are studied in
[3][27]. They have the following basic properties.

\nd{\bf Theorem (Wenzl)} Let the ground ring be a field of character
0, then $\mathcal{B}_n (\tau )$ is semisimple if and only if $\tau
\notin \mathbb{Z}$ or $\tau \in \mathbb{Z}$ and $\tau >n .$

Many algebras related to Lie theory including Hecke algebras, BMW
algebras are cellular algebras defined as follows [15]. In the same
paper Graham and Lehrer construct cellular structure for BMW
algebras and Brauer algebrsa.

\nd{\bf Definition (Graham, Lehrer)[15]} A cellular algebra over $R$
is an associative algebra $A$, together with cell datum $(\Lambda
,M,C,*)$ where
\begin{itemize}
  \item (C1) $\Lambda $ is a partially ordered set and for each $\lambda \in
\Lambda$ ,$M(\lambda)$ is a finite set such that$ C:\cap _{\lambda
\in \Lambda} M(\lambda ) \times M(\lambda ) \rightarrow A $ is an
injective map with image an R-basis of A.
 \item (C2) If $\lambda \in \Lambda $ and $S,T\in M(\lambda )$,
  write $C(S,T)= C^{\lambda } _{S,T} \in A$. Then $*$ is an
  $R$-linear anti-involution of A such that $*(C^{\lambda } _{S,T} )  =C^{\lambda }
  _{T,S}$.
\item (C3) If $\lambda \in \Lambda $ and $S,T \in M(\lambda )$
  then for any element $a\in A$ we have

 $aC^{\lambda } _{S,T} \equiv \sum _{S^{'}\in M(\lambda )} r_a (S^{'} ,S) C^{\lambda} _{S^{'},T} (modA(<\lambda))
   $

   Where $r_a (S^{'} ,S)\in R$ is independent of $T$ and where $A(<\lambda
   )$ is the R-submodule of $A$ generated by $\{ C^{\mu } _{S^{''} ,T^{''} } | \mu < \lambda ; S^{''}, T^{''} \in M(\mu )
   \}.$.

\end{itemize}

In table 1 we observe the presentation can be related to type
$A_{n-1} $ Dynkin diagram in the following way.  Every node $p_i $ (
$1\leq i\leq n-1 $) of the Dynkin diagram corresponds to a pair of
generators $s_i , e_i $. Between $s_i ,e_i $ there are relations
$s_i ^2 =1$, $e_i ^2 =\tau e_i $ and  $s_i e_i =e_i s_i =e_i $.
There are two pattern of relations between $s_i , e_i $ and $s_j
,e_j $. When $p_i $ is connected to $p_j $ by an edge, equivalently
$i=j \pm 1$, there are relations $s_i s_j s_i = s_j s_i s_j $, $e_i
s_j e_i = e_i $, $e_j s_i e_j = e_j $  and $s_i s_j e_i = e_j s_i
s_j $. When $p_i $ isn't connected to $p_j $, there are relations
$s_i s_j =s_j s_i $, $e_i e_j =e_j e_i $, $s_i e_j = e_j s_i $.
These are all relations for $B_n (\tau )$.

 A simply-laced Dynkin
diagram is a Dynkin diagram whose edges are simply-laced. Finite
type simply-laced Dynkin diagram consists ADE type Dynkin diangrams.
For every such Dynkin diagram $\Gamma$,  above discussion shows how
to define algebras similar with $B_n (\tau )$
 and $B_n (\tau ,l )$. This is due to Cohen, Gijsbers and Wales. In [10] they define  algebra $B _{\Gamma } (\tau )$ and algebra
 $B_{\Gamma } (\tau , l)$ with generators and relations as in table 2.
 These are generalization of Brauer algebra and BMW algebra to
 simply-laced type root systems (or Coxeter groups ).Let $I$ be the set of nodes of $\Gamma$. When $i, j \in I$ are
connected by an edge we write $i\sim j$. Otherwise we write $i\nsim
j$. Set $m=\frac {l-l^{-1}}{1-\tau }$.
\\

\begin{tabular}{|l|l|l|}

  % after \\: \hline or \cline{col1-col2} \cline{col3-col4} ...
\hline

   & $B_{\Gamma } (\tau )$  & $B _{\Gamma } (\tau,l)$ \\
\hline

  Generators    & $s_i  (i\in I ) $; $e_i (i\in I) $& $X_i (i\in I)$ ; $E_i (i\in I)$ \\
\hline

  Relations   &  $s_i s_j s_i =s_j s_i s_j$, $if\ i\sim j $  ;  & $X_i X_j X_{i} =X_j X_{i} X_j$, $if\ i\sim j$;
\\
 & $s_i s_j e_i = e_j s_i s_j $ if $i \sim j$; & $X_i X_j E_i = E_j X_i X_j
 $ if $i \sim j$;\\
 & $s_i s_j =s_j s_i$ ,$if\ i\nsim j $;& $X_i X_j =X_j X_i ,if\ i\nsim j $ ;\\
 &$s_i ^2 =1 ,for\  all\  i $;& $l(X_i ^2 +mX_i -1) =mE_i , for\  all\
i $;\\
 & $s_i e_i =e_i ,for\ all\ i $;&$ X_i E_i = l^{-1} E_i ,for\ all\ i $;\\
& $e_i s_{j} e_{i} = e_{i}, if\ i\sim j $; &$ E_i X_{j} E_{i} = l
E_i ,if i\sim j $;
\\

&$s_i e_j = e_j s_i , if\ i\nsim j $; & $X_i E_j =E_j X_i ,if i\nsim j $;\\
&$e_i ^2 = \tau e_i ,for\  all\  i.$ &$ E_i ^2 =\tau E_i $
. \\

  \hline
\end{tabular}
\\

These generalized Brauer algebras have no graph to representing
their elements any more, but they have almost all important
algebraic properties of Brauer algebras.

\subsection{Pseudo reflection groups, Complex braid groups and Hecke algebras}
Let $V$ be a complex linear space. An element $s$ in $GL(V)$ is
called a pseudo reflection if it can be presented as $diag(\xi
,1,\cdots ,1)$  where $\xi$ is a root of unit. If $\xi$ is $-1$ then
$s$ is simply called a reflection. $G\subset GL(V)$ is called a
pseudo reflection group if it is generated by pseudo reflections. If
$G$ is generated by reflections then it is called a complex
reflection group.  When $V$ is an irreducible representation of $G$,
$(V,G)$ is called an irreducible pseudo reflection group. Every
pseudo reflection group is isomorphic to direct product of some
irreducible factors. Isomorphism class of irreducible pseudo
reflection groups are classified by Shephard-Todd [25]. They
consists of a infinite family $\{\  G(m,p,n)\  \} $ (  $n\leq 1
$,$m\leq 2$ ,$p|m $ ) and 34 exceptional ones. Like reflection
groups, one can associate a braid group and a Hecke algebra with
every pseudo reflection group as follows.

Let $(V,G)$ be a pseudo reflection group. Let $R=\{ s_i \} _{i\in I}
$ be the set of pseudo reflections contained in $G$. For a pseudo
reflection $s$, the reflection hyperplane $H_s$ is defined as the
$1$ characteristic space of $s$, which is a codimension 1 subspace
of $V$. Let $M_G = V -\cup _{s\in R} H_s $ be the complementary
space of all reflection hyperplanes in $V$. A theorem of Steinberg
says the action of $G$ on $M_G $ is free. The braid group $A_G
$associated with $G$ is defined as $\pi _1 (M_G / G)$. It's subgroup
$P_G = \pi _1 (M_G)$ is called the pure braid group associated with
$G$.  There is a short exact sequence  $1\rightarrow P_G \rightarrow
A_G \rightarrow G \rightarrow 1 .$

As in [2][5], a Hecke algebra $H _G (\bar {\lambda})$can be
associated with every pseudo reflection group $G$, where $\bar
{\lambda}$ is a set of parameters.  The Hecke algebra $H _G (\bar
{\lambda})$ is a quotient algebra of the group algebra $\mathbb{C}
B_G $. It's dimensional equals $|G|$. For generic $\bar {\lambda}$,
$H _G (\bar {\lambda})$ is a semisimple algebra whose irreducible
representations are in one to one correspondence with those of $G$
in a natural way. This correspondence can be described by the
following Cherednik type connection.

Suppose the set of reflection hyperplanes of $G$ is $\{H_v \}_{v\in
P}$ and the set of pseudo reflections in $G$ is $R$. Define a map
$v: R\rightarrow P$ by requesting the reflection hyperplane of $s$
to be $H_{v(s)}$.  There is a natural action of $G$ on $P$ induced
by the action of $G$ on the set of reflection hyperplanes.  Now for
every reflection hyperplane $H_v$, chose a linear form $f_{v} $ with
kernel $H_v$. Define a closed 1 form
 on $M_G $ as $\omega _v = d\log f_v $.  Suppose $\{ \mu _s \} _{s\in R
 }$ is a set of constants satisfying the condition:
$\mu _{s_1 } =\mu _{s_2} $ if $s_1 $ is conjugate to $s_2$.

\begin{prop}{$ [5]$}
The formal connection  $\Omega _G  = \kappa \sum _{v\in P } (\sum
_{s\in R,v(s)=v }\mu _s s )\omega _v $ on $M_G \times \mathbb{C} G$
is  flat and $G$- invariant.

\end{prop}

Here for convenience we use a slightly different form of $\Omega _G
$. Suppose $(E,\rho )$ is a linear representation of $G$ . By above
proposition $\rho (\Omega _G ) = \kappa \sum _{v\in P } (\sum _{s\in
R,v(s)=v }) \mu _s \rho (s) \omega _v $ defines a flat connection on
the bundle $M_G \times E$. It can induce a flat connection $ \bar
{\Omega } _{\rho }$ on the quotient bundle $M_G \times _G E $
because of $G$ invariance of $\rho (\Omega _G )$. By taking
monodromy a family of representations of $B_G $ parameterized by
$\kappa $, $\mu _s $'s are obtained. It is proved in [5] that these
monodromy representations factor through $H _G (\bar {\lambda})$ for
suitable choice of $\bar {\lambda} $.

We give a simple setting for $G$-invariance of connections. Let $A$
be an algebra with unit $1$, $A^{\times }$ be the set of invertible
elements in $A$, which form a group under the multiplication. Let
$\varphi : G \rightarrow A^{\times }$ be a morphism of groups.  Let
$\rho : A\rightarrow End(E)$ be a representation of $A$ such that
$\rho (1)= id_E $. Then $\varphi \circ \rho : G \rightarrow End(E)$
is a representation of $G$, we denote it as $\bar {\rho }$.

The group $G$ acts on the bundle $E\times M_G $ as: $g\cdot (v, p) =
(\bar {\rho }(g)(v) ,g\cdot p )$, for $g\in G$,$v\in E$, and $p\in
M_G $. The quotient $E\times M_G /G$ is a linear bundle on $M_G /G$.
And $G$ acts on $A\times M_G $ as: $g\cdot (x,p) = (\varphi (g) x
\varphi (g)^{-1} ,g\cdot p )$, for $g\in G$, $x\in A$, and $p\in M_G
$. Then the quotient $A\times M_G /G $ is a bundle of algebras on
$M_G $.

Now Let $\Omega = \sum _{i\in P} X_i \omega _i $ be a connection on
$A\times M_G $, where $X_i \in A$. Then $\Omega _{\rho } = \sum
_{i\in P} \rho (X_i ) \omega _i $ is a connection on the bundle $E
\times M_G .$ We have

\begin{prop}[G-invariance]
The connection $\Omega $ induce a connection on the quotient bundle
if and only if  $X_{w(i)}= \varphi (w) X_i \varphi (w) ^{-1}$, for
any $i\in P$, $w\in G$.

The connection $\Omega _{\rho }$ induce a connection on the quotient
bundle $E\times M_G $ if and only if $\rho (X_{w(i)}) = \bar
{\rho}(w) \rho (X_i ) \bar {\rho}(w) ^{-1} $, for any $i \in P$,
$w\in G$.

\end{prop}

When $\Omega$ and $\Omega _{\rho}$ are flat connections, the
equation of flat sections for $\Omega $ is: $dF = \Omega F$, Where
$F: M_G \rightarrow A$ is a section. The equation of flat sections
for $\Omega _{\rho}$ is $d\Phi = \Omega _{\rho} \Phi $, where $\Phi
: M_G \rightarrow E$ is a section.  Flat sections of these two
bundles are connected as follows. Let $U\subset M_G $ be some
region, and $F: U\rightarrow A $ be a local flat section. Then for
any $v\in E $, $F\cdot v : U\rightarrow E$ is a local flat section
on $E\times M_G $.

We consider the monodromy operators when $\Omega $ and $ \Omega
_{\rho }$ are flat. Let $\pi : M_G \rightarrow M_G /G $ be the
quotient map. Chose a base point $p_0 \in M_G $, and chose $\bar
{p_0 }= \pi (p_0 ) $ as the base point of $M_G /G$. We define some
special path classes named generators of monodromy as follows. This
presentation is borrowed from $[5]$.

Let $H_v$ be a reflection hyperplane, $q\in H_v \setminus \sum
_{u\neq v} H_u $ be a generic point.Let $L_v =im (s_v -Id_V )$,
which is a complex line and we have $V= H_v \oplus L_v $. For $p\in
V $, we set $p= pr_{v } (p) + pr^{\perp } _{v } (p)$ with $pr_{v }
(p) \in H_v $ and $pr^{\perp } _{v } (p)\in L_v $. Thus, we have
$$ s_v (p)= e^{\frac {2\pi \sqrt{-1}}{m_v }} pr^{\perp } _{v } (p) +pr_{v }
(p).$$

 let  $\gamma : [0,1]\rightarrow V $ be a path in $V $ starting from $p_0 $ and
ending at $q $,such that $ \gamma (t) \in M_G$ for $t\neq 1$. For
$\epsilon \in (0,1]$, we denote the partial path $\gamma |
_{[0,1-\epsilon ]} : [0,1-\epsilon ]\rightarrow V $ as $\gamma
_{\epsilon}$, and the point $\gamma (1-\epsilon)$ as $p_{\epsilon}
$.

For $p\in V$, we define a path $\sigma _{v, p }$ from $p$ to $s_v
(p)$ as follows:
$$\sigma _{v,p} : [0,1]\rightarrow V,\ t\mapsto e^{\frac {2\pi \sqrt{-1}t}{m_v }} pr^{\perp } _{v } (p) +pr_{v }
(p).$$

Denote the composed path $s_v (\gamma _{\epsilon} ^{-1} )\circ
\sigma _{v, p_{\epsilon }} \circ \gamma _{\epsilon }$ as $ \sigma
_{v,\gamma ,\epsilon }$, which is a path in $V$ from $p_0 $ to $s_v
(p_0 )$. The image $\pi (\sigma _{v,\gamma ,\epsilon })$ is a loop
in $V/G $ based at $\bar {p_0 }$. It isn't difficult to see when
$\epsilon $ is small enough then the path  $ \sigma _{v,\gamma
,\epsilon }$ is in $M_G $, and $\pi (\sigma _{v,\gamma ,\epsilon })$
is in $M_G /G $. Moreover if $\epsilon $ is small enough the
homotopy class of $\pi (\sigma _{v,\gamma ,\epsilon })$ doesn't
depend on $\epsilon $, thus determine an element $s_{v,\gamma }\in
\pi _{1} (M_G /G ,\bar {p_0 })$.  We call this kind of elements in
$\pi _1 (M_G /G , \bar {p_0 })$ as generators of monodromy.

Now suppose we have flat and $G$-invariant connection $\Omega $,  We
denote the induced flat connection on the quotient bundle as $\bar
{\Omega }$. By taking monodromy we obtain an operator $ \psi :\pi _1
(M_G /G , \bar {p_0 }) \rightarrow A ^{\times }$. Let $\rho :
A\rightarrow End(E)$ be a representation of $A$. We have a flat and
$G$-invariant connection $\Omega _{\rho }$ on the bundle $E\times
M_G $, which induce a flat connection $\bar {\Omega _{\rho } }$ on
the quotient bundle $(E\times M_G )/G$. We have the monodromy
representation $\psi _{\rho } : \pi _1 (M_G , \bar {p_0 })
\rightarrow GL (E)$.  Relationship between these two kind of
monnodromy representations can is simply:
$$\rho (\psi (\lambda )) = \psi _{\rho }(\lambda ), for\ any\ \lambda \in \pi _1 (M_G , \bar {p_0
}).$$

Let $s_{v, \gamma } \in \pi _1 (M_G , \bar {p_0 })$ be a generator
of monodromy. The element $\psi (s_{v,\gamma })\in A $ and the
morphism $\psi _{\rho} (s_{v,\gamma })$ can be described as follows.
Let $\epsilon$ be small enough such that $\sigma _{v,\gamma
,\epsilon }$ lies in $M_G $ and it's homotopy class is $s_{v,\gamma
}$.  On the bundle of algebras $A\times M_G $, by parallel
transportation along the path $\sigma _{v,\gamma ,\epsilon }$, we
obtain a morphism $T _{\gamma ,\epsilon }: A_{p_0 }\rightarrow
A_{s_v (p_0 )}.$ Since when forming quotient $A_{p_0 }$ is
identified with $A_{s_v (p_0 )}$ by the correspondence $a\mapsto
\phi (s_v ) a $, so
$$\psi (s_{v,\gamma })= \phi (s_v
)^{-1} T_{\gamma ,\epsilon } (1).$$
   Similarly if $T
_{\gamma ,\epsilon , \rho } : E_{p_0 }\rightarrow E_{s_v (p_v )}$ is
parallel transportation of the connection $\Omega _{\rho }$,  we
have $$\psi _{\rho} (s_{v,\gamma }) = \rho (\phi (s_v)) ^{-1} T
_{\gamma ,\epsilon , \rho } . $$

\section{Flat connections for BMW algebras }

In this section we define a flat connection for BMW algebra $B _n
(\tau ,l)$. Using this connection we can deform any representation
of the Brauer algebra $B _n (\tau )$ to be a one parameter family of
representations of the braid group $B_n $, which will be shown to
factor through $B _n (\tau ,l)$.

We begin with some knowledge for hyperplane arrangements. Let $E$ be
a complex linear space. An hyperplane arrangement(or arrangement
simply)in $E$ simply means a finite set of hyperplanes contained in
$E$ .Arrangements arise in many fields of mathematics. For example,
let $(V,G)$ be a pseudo reflection group. Then the set of reflection
hyperplanes of $G$ is an arrangement. We denote it as $\mathscr{A}
_G $. For an arrangement $\mathscr{A} =\{ H_i \} _{i\in I}$ in $E$
let $M_{\mathscr{A}} = E- \cup _{i\in I} H_i $ be the complementary
space of all hyperplanes of $\mathscr{A}$. Topology and geometry of
$M_{\mathscr{A}}$ are the main topics of arrangement theory.

Let $\mathscr{A} =\{ H_i \} _{i\in I}$ be an arrangement in $E$. An
edge of $\mathscr{A}$ is defined to be any nonempty intersection of
elements of $\mathscr{A}$. If $L$ is an edge of $\mathscr{A}$,
define an subarrangement
$$\mathscr{A} _L = \{ H_i \in \mathscr{A} |
L\subset H_i \} =\{H_i \}_{i\in I_L } .$$

For every $i\in I$, chose a linear form $f _i $ with kernel $H_i $.
Define an 1 form on $M_{\mathscr{A}}$ as $\omega _i = d\log f_i $.
Consider the formal connection $\Omega = \kappa \sum _{i\in I} X_i
\omega _i $.  Here $X_i $ are linear operators to be determined.
When we take $X_i $'s as homeomorphisms of some linear space $E$,
then $\Omega $ is realized as a connection on the bundle
$M_{\mathscr{A}} \times E$.  We have the following theorem.

\begin{thm}{$[19]$}
 The formal connection $\Omega $ is flat if and only if the
following conditions are satisfied.

 For any codimension  2  edge  $L$  of  $\mathscr{A}$,and for any
$i\in I_L$,  $[ X_i ,\sum _{j\in I_L } X_j ] =0$ .

Here $[A,B]$ means $AB-BA$.
\end{thm}

Let $E$ be a n-dimensional Euclidean space with an orthonormal basis
$\{v_1 ,\cdots ,v_n \}$. Suppose the corresponding coordinate system
is $(x_1 ,\cdots ,x_n )$. For $1\leq i<j\leq n$, denote the
hyperplane $ker(x_i -x_j )$ as $H_{i,j}$. Let $E^{\mathbb{C}}=
E\otimes _{\mathbb{R}}  \mathbb{C} \cong \mathbb{C} <v_1 ,\cdots
,v_n >$. Denote the corresponding coordinate system as $(z_1 ,\cdots
,z_n )$. There is a natural imbedding $j: E\rightarrow
E^{\mathbb{C}}$. Denote $ker (z_i -z_j )$ as $H_{i,j }
^{\mathbb{C}}$, which is called the complexification of $H_{i,j}$.
We define  $Y_n = E\setminus \cup _{i<j} H_{i,j}$, $Y_n
^{\mathbb{C}} =E^{\mathbb{C}} \setminus \cup _{i<j} H_{i,j}
^{\mathbb{C}}$. Define $\Delta =\{ (x_1 ,\cdots ,x_n )\in E | x_1 <
x_2 <\cdots <x_n \}$ which is called a chamber. It isn't hard to see
the walls of $\Delta$ are  hyperplanes $\{H_{1,2} ,H_{2,3} ,\cdots
,H_{n-1 ,n } \}.$ Chose $p\in \Delta $. For $1\leq i\leq n-1$, chose
$p_i \in H_{i,i+1 }$ lying in the closure of $\Delta$ and not lying
in any other $H_{j,k}$. We chose path $\gamma _{i} :[0,1]
\rightarrow E$ for any $1\leq i\leq n-1 $ satisfying the following
conditions:
$$\gamma _i  (0) =p ,\ \gamma _i (1) =p_i ,\ \gamma _i (t) \in \Delta \ for\ 0<t<1 . $$

The symmetric group $S_n$ acts on $Y_n ^{\mathbb{C}}$ freely by
permuting the coordinates. Denote the quotient manifold $Y_n
^{\mathbb{C}}/ S_n $ as $X_n ^{\mathbb{C}}$. Denote the quotient map
from $Y_n ^{\mathbb{C}}$ to $X_n ^{\mathbb{C}}$ as $\pi $. We take
$j(p)$, $\bar {j(p)} = \pi (j(p))$ as base point of $Y_n
^{\mathbb{C}}$ and  $X_n ^{\mathbb{C}}$ respectively. According to
the discussion in 2.2 , the paths $p \circ \gamma _i $ can be used
to define a element $s_{\gamma _i } \in \pi _1 (X_n ^{\mathbb{C}},
\bar {j(p)})$ (generator of monodromy). The following result can be
found in [6].

\begin{prop}{$[6]$}
 The map $\sigma _i \mapsto s_{\gamma _i } $ extends to an
isomorphism $B_n \rightarrow \pi _1 (X_n ^{\mathbb{C}}, \bar
{j(p)}).$
\end{prop}

The following simple will be used later.

\begin{lem}
In the Brauer algebra $B _n (\tau )$, let $s_{i,j} \in S_n \subset B
_n (\tau )$ be $(i,j)$ permutation., let $e_{i,j}$ be as in
introduction. we have
\begin{itemize}
\item[(1)] $e_{i,j} s_{k,l} =s_{k,l} e_{i,j} $ if $ \{i,j\} \cap
\{k,l\}=\varnothing ;$
\item[(2)] $e_{i,j} e_{k,l} =e_{k,l} e_{i,j}$
if $ \{i,j\} \cap \{k,l\}=\varnothing ;$
\item[(3)] $e_{i,j} = e_{j,i};$
\item[(4)] $e_{i,j} e_{i,k} = s_{j,k} e_{i,k} = e_{i,j} s_{j,k}$ ,for
any different $i,j,k$;
\item[(5)] $ e_{i,j} ^2 = \tau e_{i,j}$ , for
any $ i\neq j$; (6) $s_{i,j} e_{j,k} = e_{i,k} s_{i,j}$.
\end{itemize}
\end{lem}

\begin{pf}
 They can be checked directly by using graghs.
\end{pf}
\\

 For $1\leq i<j\leq n-1 $, define $\omega _{i,j}= d(z_i -z_j )/(z_i -z_j
)$. Consider the formal connection $\bar{\Omega} _n = \kappa \sum
_{i<j } (s_{i,j} -e_{i,j}) \omega_{i,j}$. We have

\begin{prop}
The formal connection $\bar{\Omega} _n$ is flat and $S_n $
invariant.

\end{prop}

\begin{pf}
We certify $\bar{\Omega} _n$ satisfies conditions of theorem 3.1.
For the arrangement $\mathscr{A} _n $,there are then following two
type of codimension 2 edges

Case 1. $L= H_{i,j} \cap H_{k,l}$, $\{i,j\}\cap
\{k,l\}=\varnothing$.

Whence $\mathscr{A} _L =\{H_{i,j} ,H_{k,l}\}$. Now we have $s_{a,b}
s_{c,d} =s_{c,d} s_{a,b}$  and $e_{a,b} e_{c,d} =e_{c,d} e_{a,b}$ if
$\{a,b \}\cap \{c,d \}=\varnothing .$ They are most easily seen by
using graphs. so $[s_{i,j} -e_{i,j} ,s_{k,l} -e_{k,l}]=0.$ Which
gives
$$[s_{i,j} -e_{i,j} ,s_{i,j} -e_{i,j} +s_{k,l} -e_{k,l}]=0=
 [s_{k,l} -e_{k,l} ,s_{i,j} -e_{i,j} +s_{k,l} -e_{k,l}].$$

Case 2. $L=H_{i,j} \cap H_{j,k}$, where $i,j,k$ are different. In
this case $\mathscr{A} _L = \{ H_{i,j} ,H_{j,k} ,H_{i,k} \}$,
\begin{align*}
&[s_{i,j} - e_{i,j} ,  s_{i,k} -e_{i,k} +s_{j,k} -e_{j,k}] \\
 =&[s_{i,j},s_{i,k} +s_{j,k}] +(-e_{i,j} s_{i,k}
+e_{i,j}
e_{j,k})+(s_{i,k}e_{i,j} -e_{j,k} e_{i,j} ) \\
   +&
 (-e_{i,j} s_{j,k} +e_{i,j} e_{i,k} ) + (s_{j,k} e_{i,j}
-e_{i,k}e_{i,j}) +[s_{i,j} ,-e_{i,k}-e_{j,k}]  \\
=& (-e_{i,j} s_{i,k} +e_{i,j}e_{j,k})+(s_{i,k}e_{i,j} -e_{j,k}
e_{i,j} ) \\
 +&(-e_{i,j} s_{j,k} +e_{i,j} e_{i,k} ) + (s_{j,k}
e_{i,j}
-e_{i,k}e_{i,j})  +[s_{i,j} ,-e_{i,k}-e_{j,k}]  \\
 =& (-e_{i,j}s_{i,k} +e_{i,j}
e_{j,k})+(s_{i,k}e_{i,j} -e_{j,k} e_{i,j} )
 +(-e_{i,j} s_{j,k} +e_{i,j} e_{i,k} ) \\
  +&(s_{j,k} e_{i,j}-e_{i,k}e_{i,j})
 =0.
\end{align*}

The second equality is because $s_{i,j} s_{i,k} +s_{i,j} s_{j,k} =
s_{j,k} s_{i,j} +s_{i,k} s_{i,j}. $ For the third equality use Lemma
3.1, (6). For the fourth equality use lemma 3.1, (4). $G$-invariance
of $\bar{\Omega} _n$ is evident.

\end{pf}

Let $(E,\rho)$ be a finite dimensional representation of $B _n (m
)$. By proposition 3.1, the connection
$$ \rho (\bar{\Omega} _n)= \kappa \sum _{i<j }
( \rho (s_{i,j} ) -\rho (e_{i,j}) ) \omega_{i,j}$$ induce a flat
connection on the bundle $Y_n \times _{S_n} E $, which is a linear
bundle on $X_n$. So by taking monodromy we obtain representations of
the braid group $B_n$. We have

\begin{thm}
The monodromy representations of $B_n$ constructed above factor
through $B _n (\tau , l)$, for $\tau = \frac{q^{1-m} - q^{m-1} +
q^{-1} -q}{q^{-1} -q}$, $l=q^{m-1} $. Where $q= \exp \kappa \pi
\sqrt{-1}$.
\end{thm}

\begin{pf}
We first prove this theorem in case of $n=2,3.$ The cases for $n\geq
4$ can be reduced to cases for $n=2,3$ by a result in local theory
of meromorphic connections.

Case 1.  $(n=2)$  The Brauer algebra $B _2 (m )$ is 3 dimensional,
with a basis $\{ 1,s,e \}$, submit to relations $e^2 =m e , se=es=e
,s^2 =1 .$ where $1$ means the unit. Idempotents of this algebra are
$$\varepsilon _0 = \frac {e}{m} ,\ \varepsilon _1 = \frac {1-s}{2}
,\ \varepsilon _2 = \frac {1+s}{2} -\frac {e}{m} .$$

In this case the relevant reflection group is the symmetric group
$(S_2 ,  \mathbb{C} ) $, the only reflection hyperplane is $0\in
\mathbb{C}$. The complementary space is $\mathbb{C} ^{\times }$. Let
$p_0 =1\in \mathbb{C} ^{\times } $ be a base point. Define
$$\sigma : [0,1]\rightarrow \mathbb{C} ^{\times } , t\mapsto e^{t\pi \sqrt{-1 }}
$$
which is a path from $1$ to $-1= s(1)$. Let $\pi :\mathbb{C}^{\times
} \rightarrow \mathbb{C}^{\times} /S_2 $ be the quotient map. Then
the loop $\pi (\sigma )$ represent a generator of $\pi _1
(\mathbb{C}^{\times} /S_2 ,\bar {p_0 })\cong \mathbb{Z},$ which is
denoted as $s_{\sigma }$.

The flat connection is $\Omega = \kappa (s-e) dz/z .$ Let $F$ be a
section of the bundle of algebras $\mathbb{C} ^{\times } \times B_2
(m)$, which can be presented as $F= F_0 \varepsilon _0 +F_1
\varepsilon _1 + F_2 \varepsilon _2 .$ Where $F_i$'s are functions
on $\mathbb{C} ^{\times }$. $F$ is a flat connection when
 $dF = \Omega F .$ It can be written as:
$$ \frac {dF_0}{dz} = \frac {\kappa (1-m)}{z} F_0 ,\  \frac {dF_1}{dz}
= \frac {-\kappa }{z} F_1 ,\  \frac {dF_2}{dz} =  \frac {\kappa }{z}
F_2 . $$

The solution of this equation is $F(z) = z^{\kappa (1-m )} c_0
\varepsilon _o +z^{-\kappa} c_1 \varepsilon _1 + z^{\kappa} c_2
\varepsilon _2 .$ where $c_i $ are constants. Chose $1\in \mathbb{C}
^{\times }$ as the base point. If we request the initial condition
$F(1) = 1$, then $c_0 =c_1 =c_2 =1 .$ Through parallel
transportation along the path $\sigma $, $1$ is corresponded to the
element $T_{\sigma }(1)=e^{\kappa (1-m )\pi \sqrt{-1 }} \varepsilon
_o +e^{-\kappa\pi \sqrt{-1 }}  \varepsilon _1 + e^{\kappa\pi
\sqrt{-1 }}  \varepsilon _2 $.

Let $\psi : \pi _1 (\mathbb{C}^{\times} /S_2 ,\bar {p_0
})\rightarrow B_2 (m ) ^{\times }$ be the monodromy morphism. We
have
$$\psi (s_{\sigma })= s^{-1} T_{\sigma }(1)  = e^{\kappa (1-m )\pi \sqrt{-1 }}  \varepsilon
_o -e^{-\kappa\pi \sqrt{-1 }}  \varepsilon _1 + e^{\kappa\pi
\sqrt{-1 }}  \varepsilon _2 .$$

Since $\varepsilon _i $'s are idempotents, and $\sum _{i=0} ^{2}
\varepsilon _i =1$ and $\varepsilon _i \varepsilon _j =0$ for $i\neq
j$, so we have
$$(\psi (s_{\sigma })-e^{\kappa (1-m )\pi \sqrt{-1 }})(\psi (s_{\sigma })+e^{-\kappa\pi \sqrt{-1 }})(\psi (s_{\sigma })-e^{\kappa\pi
\sqrt{-1 }}) =0.$$

By definition the BMW algebra $B_2 (\tau ,l)$ is 3 dimensional with
a basis $\{1, X,E \}$ submitting to the relations
$$XE= l^{-1} E =EX ,\ E^2 =\tau E,\  l(X^{2} +JX -1) =JE.$$
where $J= \frac {l-l^{-1}}{1-\tau }.$

Let $q= e^{\kappa \pi \sqrt{-1}}$, $l =e^{\kappa (m -1) \pi
\sqrt{-1}}$ and $\alpha = q^{-1}-q $, $\tau =\frac{ l^{-1} -l
+\alpha }{\alpha} = \frac{q^{1-m} -q^{m-1}+q^{-1} -q}{q^{-1} -q}$
From theorem 4.1, we see $\psi $ factor through $\iota : \pi _1
(\mathbb{C}^{\times} /S_2 ,\bar {p_0 })\rightarrow B_2 (\tau ,l )$
because if we set
$$X= \psi (s_{\sigma }),\ E= \frac {l }{q^{-1} -q}
(\psi (s_{\sigma })^2 + (q^{-1} -q)\psi (s_{\sigma }) -1),$$ then
all relations for BMW algebras are satisfied.

Case 2. $n=3$. We only consider those $B_3 (m)$'s which are
semisimple. We prove for any finite dimensional representation $\rho
:B_3 (m)\rightarrow End(E )$, the monodromy representation of $\rho
(\bar{\Omega} _3): B_3 \rightarrow Aut (E)$ factor through $\iota :
B_3 \rightarrow B_3 (\tau ,l )$, where $\tau $, $l $ are defined as
in case 1. For this aim we only need to consider the cases when
$\rho$ is an irreducible representation. It isn't hard to see that
$B_3 (m) $ has 4 indifferent irreducible representations as follows:

$\rho _1 $: $E= \mathbb{C} <v >$. $s_i v =v $ for any $i$; $e_i v=0
$ for any $i$.

$\rho _2 $: $E= \mathbb{C} <v> $. $s_i v=-v $ for any $i$; $e_i v=0
$ for any $i$.

$\rho _3 $: $E= \mathbb{C} <v_1 ,v_2 >$. $e_i v_j =0 $ for any
$i,j$. The actions of $s_i $'s make $E$ becoming the canonical
representation of $S_3 $.

$\rho _4 $: $E= \mathbb{C} <v_{1,2} ,v_{2,3} ,v_{1,3} >$. $s_{i,j}
v_{i,j} =v_{i,j}$. For $i>j$, we set $v_{i,j}=v_{j,i}$ and
$s_{i,j}=s_{j,i}$. For different $i,j,k$, $s_{i,j} v_{j,k} =v_{i,k
}$. $e_{i,j} v_{i,j}= m v_{i,j}$. For different $i,j,k$, $e_{i,j}
v_{j,k}= v_{i,j}.$

In the case of $\rho _i $, let $\psi _i : \pi _1 (X_n ^{\mathbb{C}},
\bar {j(p)})\rightarrow Aut (E)$ be the monodromy representation of
the connection $\rho _i (\bar{\Omega} _3 )$.  For $k=1,2$, we set
$X_k = \psi _i (s_{\gamma_k })$ and $E_k = \frac {l }{q^{-1} -q}
(\psi (s_{\sigma _k })^2 + (q^{-1} -q)\psi (s_{\sigma _k }) -1)$.

In the case of $\rho _1$, $\rho _2 $, $\rho _3$, the elements $e_i
$'s act as 0, so the connection $\rho _i (\bar{\Omega} _3 )$ in fact
degenerates to the flat connection for Hecke algebras. It is readily
to check in these cases $\psi (s_{\sigma _k })^2 + (q^{-1} -q)\psi
(s_{\sigma _k }) -1 =0$. So $E_k =0 $ in these cases.

In the case of $\rho _4 $, we notice this is nothing but the
infinisimal Krammer representation defined by Marin.  This case is
proved by using theorem 4.1 in $[21]$.

\end{pf}

\section{Generalized Lawrence-Krammer Representations }

Let $V$ be a $n$-dimensional complex linear space. Let $G\subset
GL(V)$ be a finite reflection group. Let $R$ be the set of
reflections in $G$. For each $s\in R$, denote the reflection
hyperplane of $s$,namely the subspace $Ker(s-1)$ as $H_s$. Let $M_G
= V\backslash \bigcup _{s\in R} H_s $ be the complementary space of
all reflection hyperplanes.

The generalized LK representations of $B_G$ of Marin are described
by certain flat connections as follows.  First, for every $H_s $,
let $\alpha _s$ be a nonzero linear form with kernel $H_s $. Let
$\omega _s = d \alpha _s /\alpha _s $, which is a closed,holomorphic
1-form on $M_G $. Then let $V_G =\mathbb{C}\langle v_s \rangle
_{s\in R}$ be a complex linear space with a basis indexed by $R$.
For every pair of elements $s,u \in R$, define a nonnegative integer
$\alpha (s,u) =\# \{ r\in R | rur=s \} $. Chose a constant $m\in
\mathbb{C}$. For any $s\in R$, define a linear operator $t_s \in
GL(V_G )$ as follows:
$$t_s \cdot v_s = m v_s , \ \ t_s \cdot v_u = v_{sus} -\alpha (s,u)v_s \ \ for\ \  s\neq u. $$

Chose another constant $k\in \mathbb{C}$. Define a connection
$\Omega _{K} = \sum _{s\in R} k\cdot t_s \omega _s $ on the trivial
bundle $V_G \times M_G $.

\textbf{Theorem and Definition }([22]) The connection $\Omega _{K} $
is flat and $G$-invariant. So it induce a flat connection
$\bar{\Omega }_{K}$ on the quotient bundle $(V_G \times M_G )/G$.
The generalized LK representation for $B_G $ is defined as the
monodromy representation of $\bar {\Omega }_{K}$.

We denote the generalized Krammer representation as $(V,\rho
_{\kappa ,m })$. When $(G,V)$ is the reflection group $W_{\Gamma }$
of ADE type,  they were first constructed in $[9]$ by Cohen, Wales
and by Digne in $[13]$.  They are proved to factor through BMW
algebras in $[10]$.

\begin{thm}{ $[21] $ } The generalized Krammer representation $(V, \rho _{\kappa ,m })$factor
through the generalized BMW algebra $B _{\Gamma } (\tau ,l)$ with
$\tau =\frac {q^m -q^{-m} +q^{-1} -q}{q^{-1}-q}$ and $l= q^{-m}$.
Where $q=e^{\kappa \pi \sqrt{-1}}$.

\end{thm}

For later convenience we change notations slightly. For $s\in
\mathcal {R}$, we define $p_s : V_G \rightarrow V_G $ by
$$ p_s (v_s )= (1-m_s )v_s ,\ p_s (v_u )= \alpha (s,u ) v_s \ for\ u\neq s.$$

We also define $\iota :G \rightarrow Aut (V_G )$ by $\iota (w) (v_s)
= v_{wuw^{-1}}$.  Then $p_s$ is a projector to the complex line
$\mathbb{C}v_s $. And Marin's flat connection $\Omega _K$ is written
as $\sum _{s\in \mathcal {R}} k\cdot (\iota (s) -p_s ) \omega _s$.

Now we let $V$ be a $n$-dimensional complex linear space. Now let
$G\subset GL(V)$ be a finite pseudo reflection group(not only
reflection group). Let $R$ be the set of pseudo reflections in $G$.
For $s\in R$, let $H_s $ be the hyperplane fixed by $s$. Let
$\mathscr{A} =\{ H_i \}_{i\in P}$ be the arrangement of reflection
hyperplanes of $G$. Define $V_G =\mathbb{C} <v_i
>_{i\in P}$. Since $w(H_v )$ is annother reflection hyperplane for any
$w\in G $ and $v\in P$, there is an action of $G$ on $P$ which
induce a representation $\iota : G\rightarrow Aut (V_G )$. We also
define $w(i)$ by $H_{w(i)}= w(H_v )$.

For $i\in P$, let $p_i :V_G \rightarrow V_G $ be a projector to
$\mathbb{C} v_i $ which is written as:
$$p_i (v_i )=m_i v_i ,\ p_i (v_j ) =\alpha _{i,j} v_i .$$

As in Section 2.2 let $\{\mu _s \} _{s\in R }$ be a set of constants
such that: $\mu _{s_1} =\mu _{s_2 } $ if $ s_1$ is conjugate to $s_2
$ in $G$. Define a function $v: R\rightarrow P $ such that
$H_{v(s)}$ is the reflection hyperplane of $s$ for any $s\in R$.

Consider a connection $\Omega _{LK} $ on the trivial bundle $V_G
\times M_G $ which have the form
$$\kappa \sum _{v\in P} (\sum _{s: v(s)=v } \mu _s \iota (s) -p_v ) \omega _v .$$

 The following theorem generalize Marin's construction to psuedo
reflection groups and with more parameters in some cases.

\begin{thm}  The connection $\Omega _{LK}$ is flat and $G$ invariant if
and only if the
 the following conditions are satisfied:\\
 $(1)$ $m_i =m_j$ if there is $w\in G$ such that $\iota (w)(v_i)=v_j$.  \\
 $(2)$ $\alpha _{i,j} =\sum _{s : \iota (s )(v_j )=v_i } \mu _{s} .$

 \end{thm}

 \begin{proof} Suppose $\Omega _{LK}$ is a flat, $G$ invariant connection.
 Let $w\in G$,
 $$w^* (\Omega _{LK})=\sum _{i\in I } (\sum _{s: i(s) =i } \mu _s \iota (w) \iota (s) \iota (w)^{-1} +
 \iota _{w} p_i \iota _{w} ^{-1} ) \omega _{w(i)} .$$

 The connection $\Omega _{LK}$ being $G$ invariant means: for any $w\in G$,
 $w^* (\Omega _{LK}) = \Omega _{LK}.$ Because $\{ \omega _i \} _{i\in I }$
 constitute a basis of $H^1 _{DR} (M_G , \mathbb{C})$, so it is
 equivalent to $\iota _{w} p_i \iota _{w} ^{-1} = p_{w(i)}$, which is
 equivalent to: $m_i = m_{w(i)} .$

Let $L$ be any codimension 2 edge of the arrangement $\mathscr{A}$.
Let $H_{i_1},\cdots , H_{i_{N}}$ be all the hyperplanes in
$\mathscr{A}$ containing $L$. The condition of $\Omega$ being flat
is written as :
\begin{equation*}
[ \sum _{s: i(s)=i_a } \mu _s \rho (s) -p_{i_a } , \sum _{v=1} ^{N}
( \sum _{s:i(s) =i_v } \mu _s \rho (s) -p_{i_v } ) ] =0.
\end{equation*}

for $1\leq a\leq N $. We discuss the case $a=1$ first. Since $\Omega
_0 $ is a flat connection, it is equivalent to:

\begin{equation*}
[ p_{i_1 } , \sum _{v=1} ^{N} ( \sum _{s:i(s) =i_v } \mu _s \iota
(s) p_{i_v } ) ] + [ \sum _{s: i(s)=i_1 } \mu _s \iota (s)  , \sum
_{v=1} ^{N} p_{i_v } ) ] =0.
\end{equation*}

Now for those $s$ such that $ i(s) =i_1$ we have $\{ s(i_1 ), \cdots
, s(i_N ) \} =\{ i_1 ,\cdots ,i_N \}. $  So we have:
\begin{equation}
\begin{split}
 [ \sum _{s: i(s)=i_1 } \mu _s \iota (s)  ,
\sum _{v=1} ^{N} p_{i_v } ) ] &=
\sum _{s: i(s)=i } \mu _s \sum _{v=1} ^{N} (\iota (s) p_{i_v }-p_{i_v }\iota (s) )\\
   &= \sum _{s: i(s)=i } \mu _s \sum _{v=1} ^{N} (  p_{s(i_v )} \iota (s) -p_{i_v }\iota (s)) \\
   &= 0. \\
\end{split}
\end{equation}
So the identity (2) is equivalent to:

\begin{equation}
\begin{split}
  [ p_{i_1 } , \sum _{v=1} ^{N} ( \sum _{s:i(s) =i_v } \mu _s \iota (s)
+p_{i_v } ) ] &= \sum _{v=2} ^{N} \sum _{s:i(s) =i_v } [ p_{i_1 } ,
( \sum _{s:i(s) =i_v } \mu _s \iota (s) +p_{i_v } ) ]\\
&=0.
\end{split}
\end{equation}

This is because $[p_{i_1 } , s]=0 $ if $i(s) =i_1 $. After splitting
the Lie bracket in equation (2), the sum of all those terms mapping
to $\mathbb{C} v_{i_u }$ is $p_{i_u } p_{i_1} -\sum _{s: s(v_{i_1 }
) =v_{i_u }}\mu_s \iota (s) p_{i_1 } $. It must be 0.   Chose $s_0 $
such that $s_0 (v_{i_1 }) =v_{i_u }$, then e have $ p_{i_u } p_{i_1
}= \alpha _{i_u ,i_1 } \iota (s_0 ) p_{i_1 }$. More over, for any
$s$
 such that $s(v_{i_1 } ) =v_{i_u }$  we have $\iota (s) p_{i_1 } = \iota (s_0 )
p_{i_1 }$. Put these identities in equation (2), we get
$$(\alpha _{i_u ,i_1 } -\sum _{s: s(v_{i_1 })=v_{i_u }} \mu _s ) \iota (s_0 ) p_{i_1 } =0 .$$

So we have $\alpha _{i_u ,i_1 } =\sum _{s: s(v_{i_1 })=v_{i_u }} \mu
_s .$

Now suppose conditions is satisfied, by the same arguments we only
need to prove above equation (4) to show $\Omega $ is flat. The
conditions (2) implies
\begin{equation}
p_i p_j = \sum _{s: s(j)=i } \mu _s \iota (s) p_j , for\ any\ i\neq
j.
\end{equation}
It also implies
\begin{equation}
p_i p_j = \sum _{s: s(j)=i } \mu _s p_i \iota (s)  , for\ any\ i\neq
j.
\end{equation}
since $\iota (s) p_j = \iota (s)p_j \iota (s) ^{-1} \iota (s) = p_i
\iota (s)$ for those $s$ such that $s(j)=i$.  Now the right hand
side of equation (4) can be written as
$$\sum _{v=2} ^N (p_{i_1 } p_{i_v } - \sum _{s: s(i_v ) =i_1 } \mu _s p_{i_1 } \iota (s))+
\sum _{v=2} ^{N} (p_{i_v } p_{i_1 } -\sum _{s: s(i_v )=i_1 } \mu _s
\iota (s) p_{i_v }).$$

So the equation (2) is true by using identities (3),(4). The $G$
invariance of the connection isn't hard to see.

\end{proof}

\textbf{Remark} In the connection $\Omega _{LK}$  if make $\mu _s
=\kappa $ for all $s$ and $m_i =m$ for all $i$ then we obtain
Marin's connection. Above theorem give another proof that the
connection is flat. It explains why the numbers $\alpha _{i,j}$ will
appear in Marin's construction, also show the LK representation is
unique in some sense. We can define and describe the generalized LK
representation as follows.

\begin{defi}[Generalized LK representation for all complex braid groups]
Following notations introduced above. Since $\Omega _{LK}$ is
$G$-invariant, it induces a flat connection $\bar{\Omega } _{LK}$ on
the quotient bundle $V_G \times _{G} M_G $.  the generalized Krammer
representation of the braid group $A_G $ is defined as the monodromy
representation of $\bar{\Omega } _{LK} $.

\end{defi}

\section{Generalization of Brauer Algebras }

\subsection{Definition and Main Properties }
 Let $(V,G)$ be any finite pseudo reflection group. Denote the set of
pseudo reflections in $G$ as $R$, and let $\mathcal {A} = \{H_i \}
_{i\in P}$ be the set of reflection hyperplanes. For $i,j \in P$,
Let $R(i ,j ) =\{s\in R\ |\  s(H_j ) =H_i  \ \}$. $For i\in P$, as
in 2.2 let $G_i $ be the subgroup of $G$ consisting of elements that
fixing all points in $H_i $, let $m_i = |G_i |$, and let $s_i $ be
the unique element in $G_i $ with exceptional eigenvalue $e^{\frac
{2\pi \sqrt{-1}}{m_i }}$. we write $s_1 \sim s_2 $ for $s_1 ,s_2 \in
R$ if they are in the same conjugacy class, and write $i\sim j $ for
$i,j \in P$ if $w(i) =j $ for some $w\in G$. Chose $\mu _s \in
\mathbb{C}$ for every $s\in R$ and $m_i \in \mathbb{C}$ for every
$i\in P$ such that
$$\mu _{s_1 } =\mu_{s_2 } \ if\ s_1 \sim s_2 ,
 m_i =m_j \ if\ i\sim j.$$

The data $\{\mu _s ,m_i \}_{s\in \mathcal {R} ,i\in P }$ will be
denoted by one symbol $\Upsilon $. We define an algebra $B _G
(\Upsilon )$ as follows.

\begin{defi}
 The algebra $B_G (\Upsilon )$ associated with pseudo reflection group $(V,G)$ is generated by
  $\{ \bar {w} \}_{w\in G } \cup  \{ e_i \}_{i\in P}$, submitting to the following
relations.
\begin{itemize}
\item[(0)] $\bar {w_1} \bar {w_2} = \bar {w_3 }$ if $w_1 w_2 =w_3 $;
\item[(1)] $ \bar{s }_i  e_i  = e_i \bar{s}_i = e_i $, for $i\in P$;
$s_i \in G$ is any pseudo reflection with reflection hyperplane
$H_i$.

\item[(2)] $e_i ^2 = m_i e_i $ ;
\item[(3)] $\bar {w} e_j  = e_i \bar {w} $ , if $w \in G  $ satisfies $w(H_j )=H_i$;
\item[(4)] $e_i e_j = e_j e_i $,
 if  $\{ \ k\in P |\  H_k \supset H_i \cap H_j \}
= \{i,j \}$ ;
\item[(5)] $ e_i e_j = (\sum _{s\in R(i,j) } \mu _s s )e_j $ ,
 if $\{ k\in P | H_k \supset H_i \cap H_j \}   \neq
  \{i,j \},$
  and $R(i,j) \neq \emptyset $.
\item[(6)] $e_i e_j =0 $, if $\{  k\in I | H_k \supset H_i \cap H_j
\} \neq  \{i,j \} $ and $R(i,j) =\emptyset $.
\end{itemize}

\end{defi}

\begin{rem}
Relation $(0)$ is nothing but letting $\mathbb{C} G$ be imbedded in
$B_G (\Upsilon )$. Relation $(3)$ along with relation $(5)$ implies
$ e_i e_j = e_i (\sum _{s\in R(i,j) } \mu _s s).$

When $(G, V) $ is a complex reflection group, there is a bijection
from $R$ to $\mathcal {A} $: $s \mapsto H_s $. So we can use $R$ as
the indices set of reflection hyperplanes, and $\mathcal {A} =\{H_s
\} _{s\in R }$. In these cases, for $s_1 ,s_2 \in R $, $R(s_1 ,s_2
)= \{ s\in R | s(H_{s_2}) =H_{s_1 } \} =\{s\in R | ss_2 s=s_1 \} $.

When $G$ is a infinite Coxeter group, by using geometric
representation of $G$ above definition still make sense. $B_G
(\Upsilon)$ for those cases can also be defined by a canonical
presentation as shown in section 5.3.

\end{rem}

First we have

\begin{thm}

  $B _G
(\Upsilon )$  is a finite dimensional algebra. And $w \mapsto \bar
{w} $ for $w\in G $ induce an injection  $j: \mathbb{C} G
\rightarrow B _G(\Upsilon )$.

\end{thm}

\begin{proof}
First by using relation (3), we can identify any word made from the
set $\{w\in G \} \coprod \{e_i \} _{i\in P } $ to a word of the form
$we_{i_1} e_{i_2 }\cdots e_{i_k }$ where $w\in G $. We call $e_{i_1}
e_{i_2 }\cdots e_{i_k }$ as the '$e$-tail' of the word $we_{i_1}
e_{i_2 }\cdots e_{i_k } $. Then if two neighboring $e_{i_v }$ and
$e_{i_{v+1 }}$ don't commute with each other, then conditions in (5)
are satisfied as can be seen in the next lemma.

\begin{lem}
If two pseudo reflection $s_1 $ and $s_2 $ don't commute with each
other, suppose the reflection hyperplane of $s_1 (s_2 )$ is $H_{i_1
} (H_{i_2 } ) $, then $\{i_1 ,i_2  \} \subsetneq \{k\in P | H_k
\supseteq H_{i_1 } \cap H_{i_2 } \} $.
\end{lem}

\begin{proof} We suppose $\{i_1 ,i_2  \} = \{k\in P | H_k \supseteq H_{i_1 }
\cap H_{i_2 } \} $. Let $L= H_{i_1 }\cap H_{i_2 }$, and $<\ ,\  >$
being a $G$-invariant inner product on $V$. Chose $v_k \in H_{i_k }$
such that $v_k \perp L $ according to $<\ ,\ >$ for $k=1,2$. Suppose
$\{ v_3 ,\cdots , v_N \}$ is a basis of $L$, then $\{ v_1 ,v_2
,\cdots ,v_N \}$ is a basis of $V$. Now since $s_1 (H_{i_2 })$ is
another reflection hyperplane containing $L$ and $s_1 (H_{i_2 })\neq
H_{i_1 }$, so we have $s_1 (H_{i_2 }) =H_{i_2 }$, which implies
$s_1$ can be presented as a diagonal matrix according to the basis
$\{v_1 ,\cdots ,v_N \}$. Similarly $s_2 $ can be presented by a
diagonal matrix according to the same basis. So $s_1 s_2 =s_2 s_1 $
which is a contradiction.
\end{proof}

The first statement follows from the next lemma.

\begin{lem}
The algebra $B _G (\Upsilon )$ is spanned by the set
$$\{w\in G_m \} \coprod \{w e_{i_1 }\cdots e_{i_ M}
| w\in G ,\ e_{i_u } e_{i_v } =e_{i_v } e_{i_u } ; \ i_v \neq i_u \
if\ u\neq v \ ; M\geq 1 \}$$

\end{lem}

\begin{proof} Let $A $ be the space in $B _G (\Upsilon )$ spanned
by elements listed in the lemma. By applying (2) in definition 5.1,
we see the algebra is spanned by elements having the form $x=we_{i_1
}\cdots e_{i_K }$ where $w\in G_m $. For convenience we call $e_{i_1
}\cdots e_{i_K }$ as the e-tail of the word $x$, and $K$ as the
length of it's e-tail. We prove that every such element belongs to
$A$ by induction on the length of their e-tails $K$. First this is
true if $K=1$. Suppose it is true for $K\leq M $. Now suppose
$x=we_{i_1 }\cdots e_{i_K }$ such that $K=M+1$.  If there are two
neighboring $e_{i_v }$, $e_{i_{v+1 }}$ don't commute, then lemma 5.1
enable us to apply (5) or (6) in definition 5.1 to identify $x$ with
a linear sum of words whose e-tail length are smaller than $M+1$.
Suppose all $e_{i_v }$'s in $x$ commute with each other, if there
are $v_1 ,v_2 $ such that $i_{v_1 } =i_{v_2 } $, we use
transpositions between $e_{i_v }$'s to
 identify $x$ with a word $y=we_{j_1 }\cdots e_{j_{M+1 } }$ such that
$j_1 =j_2 $. So $x=y=m_{j_1 }w e_{j_2 }\cdots e_{j_{M+1 } } $. If
all $e_{i_v }$'s commute and all $i_v $'s are different then $x\in A
$, and induction is completed.

 For the second statement, it isn't
hard to see the following map
$$s \mapsto s,\ for\ s\in R;\ e_i \mapsto\ 0, \ for \ i\in P$$
extends to a surjection $\pi : B _G (\Upsilon ) \rightarrow
\mathbb{C} G $, and $\pi \circ j =id $. So $j$ is injective.
\end{proof}

This completes the proof of Theorem 5.1.
\end{proof}

By Theorem 5.1, $\mathbb{C} G$ is naturally embedded in $B_G
(\Upsilon )$. From now on we will denote $\bar {w} $ in $B_G
(\Upsilon )$ simply as $w$.

The next lemma reduce one parameter in $B _G (\Upsilon )$.

\begin{lem}
For $\lambda \in \mathbb{C} ^{\times }$, Let $\mu ^{\prime } _s
=\lambda \mu _s $ for $s \in R$, and Let $m^{\prime } _i =\lambda
m_i $ for $i\in P$. Let $\Upsilon ^{\prime }=\{\mu ^{\prime } _s ,
m^{\prime } _i \}_{s \in R ,i\in P}$, then $B _G (\Upsilon ^{\prime
})\cong B _G (\Upsilon )$.

\end{lem}

\begin{proof}
Denote the generators of  $B _G (\Upsilon ^{\prime })$ appeared in
Definition 5.1 as $S_{\alpha } $'s and $E_i $'s . Then
$$S_{\alpha }\mapsto s_{\alpha } , \ E_i \mapsto e_i \ for \ \alpha \in \Phi ,\ and\ i\in P $$
extend to an isomorphism from $B _G (\Upsilon ^{\prime })$ to $B _G
(\Upsilon )$.

\end{proof}

The following lemma can be found in [22].
\begin{lem}
For two different hyperplane $H_i ,H_j \in \mathcal {A}$, If $s\in R
$ satisfies $s(H_j ) =H_i $, then $s$ fix all points in $H_i \cap
H_j $.

\end{lem}

\begin{proof}
Let $ <\ ,\ > $ be a $G$ invariant, positive definite Hermitian form
on $V$.  Let $\epsilon $ be the exceptional eigenvalue of $s$, and
let $\alpha $ be an eigenvalue of $s$ with eigenvalue $\epsilon$.
Let $\alpha _i $, $ \alpha _j $ be some nonzero vector perpendicular
to $H_i $, $H_j $ respectively. Then $\alpha \perp H_s $. The action
of $s$ can be written as $s(v) = v- (1-\epsilon ) \frac {<v, \alpha
>}{<\alpha ,\alpha >} \alpha $.  Now $s(H_j )=H_i $ implies
$s(\alpha _j ) = \alpha _j - (1-\epsilon ) \frac {<\alpha _j ,
\alpha >}{<\alpha ,\alpha >} \alpha  =\lambda \alpha _i $ for some
$\lambda \neq 0$. Denote $(1-\epsilon ) \frac {<\alpha _j , \alpha
>}{<\alpha ,\alpha >}$ as $\kappa $. The condition that  $H_i $ is
different from $H_j $ implies $\kappa \neq 0$. So we have $\alpha =
\frac {1}{\kappa} (\alpha _j - \lambda \alpha _i ) \perp H_i \cap
H_j $, and $s$ fix all points in $H_i \cap H_j $.

\end{proof}

There is a natural anti-involution  on $B_G (\Upsilon )$ which may
be used to construct a cellular structure.

\begin{lem}

The correspondence $$ w\mapsto w^{-1}\   for\  w\in G \subset B_G
(\Upsilon ),\  e_i \mapsto e_i \  for\  all\  i\in P $$ extends to
an anti-involution $*$ of $B_G (\Upsilon )$ if $\mu _s = \mu _{s^{-1
}}$ for any $s\in R $.

\end{lem}

\begin{proof}
We only need to certify $*$ keep all relations in definition 5.1. As
an example for relation (5), on the one hand $* (e_i e_j ) = e_j e_i
$, on the other hand  $* [ (\sum _{s\in R(i,j)} \mu _s s ) e_j ] =
e_j (\sum _{s\in R(i,j)} \mu _s s^{-1 } ) = (\sum _{s\in R(i,j)} \mu
_s s^{-1}) e_i = (\sum _{s\in R(j,i) } \mu _{s^{-1} } s ) e_i =
(\sum _{s\in R(j,i) } \mu _{s } s ) e_i .$
\end{proof}

Next we construct some connections for $B _G (\Upsilon )$. Suppose
$\rho : B _G (\Upsilon ) \rightarrow End(E)$ is a finite dimensional
representation. On the vector bundle $E\times M_G $, we define a
connection
$$ \Omega  _{\rho } = \kappa \sum _{i\in P } (\sum _{s:i(s)=i }\mu _s \rho (s) -\rho (e_i )) \omega _i$$

where $\kappa \in \mathbb{C}$. Let $G$ acts on $E\times M_G $ as
$w\cdot (v,x)= (\rho (w) v ,wx)$ for $w\in G$ and $(v,x)\in E\times
M_G$. Then we have

\begin{prop}
The connection $\Omega  _{\rho } $ is flat and $G$-invariant.

\end{prop}

\begin{proof}
 It is enough to deal with the case $\kappa =1$. The connection
$\Omega _{\rho } $ is $G$-invariant if and only if for any $w\in G$,
\begin{equation}
\sum _{i\in P } (\sum _{s:i(s)=i }\mu _s \rho (w) \rho (s) \rho
(w)^{-1} - \rho (w)\rho (e_i ) \rho (w)^{-1}) \omega _{w(i)} = \sum
_{i\in P } (\sum _{s:i(s)=i }\mu _s \rho (s) -\rho (e_i )) \omega _i
.
\end{equation}

By definition 5.1 (3), we have $\rho (w)\rho (e_i ) \rho
(w)^{-1}=\rho (we_i w^{-1})=\rho (e_{w(i)})$. And $\{wsw^{-1} | i(s)
= i \}= \{s| i(s)= w(i) \}$, so identity (7) follows.

Let $L$ be any codimension 2 edge for the arrangement $\mathcal
{A}$, and let $H_{i_1 }, \cdots , H_{i_N } $ be all the hyperplanes
in $\mathcal {A}$ containing $L$. By theorem 3.1, to prove $\Omega
_{\rho }$ is flat we need to show for any $u$
\begin{equation}
[ \sum _{s: i(s)=i_u } \mu _s \rho (s) -\rho ( e_{i_u } ) , \sum
_{v=1} ^{N} ( \sum _{s:i(s) =i_v } \mu _s \rho (s) -\rho (e_{i_v }
)) ] =0.
\end{equation}

Now remember the connection $ \kappa \sum _{i\in P } (\sum
_{s:i(s)=i }\mu _s \rho (s)  \omega _i )$ is flat by proposition
2.1, so (8) is equivalent to

\begin{equation}
-[ \sum _{s: i(s)=i_u } \mu _s \rho (s) , \sum _{v=1} ^{N} \rho
(e_{i_v } ) ]- [ \rho ( e_{i_u } ) , \sum _{v=1} ^{N} ( \sum
_{s:i(s) =i_v } \mu _s \rho (s) -\rho (e_{i_v } )) ]=0
\end{equation}

Because for any $s$ such that $i(s)= i_u $, there is $\{ s(H_{i_1
}), \cdots ,s(H_{i_n })\}= \{H_{i_1 }, \cdots ,H_{i_N } \}$, so
\begin{equation}
\begin{split}
\rho (s) \sum_{v=1 } ^N \rho (e_{iv }) - \sum_{v=1 } ^N \rho (e_{iv
})\rho (s) &=(\rho (s) \sum_{v=1 } ^N \rho (e_{iv }) \rho (s)^{-1} -
\sum_{v=1 } ^N \rho (e_{iv }) )\rho (s)\\
&
= (\sum_{v=1 }^N \rho (e_{s(i_v) }) - \sum _{i=1} ^N e_{i_v })=0
\end{split}
\end{equation}
So (9) is equivalent to
\begin{equation}
[ \rho ( e_{i_u } ) , \sum _{v=1} ^{N} ( \sum _{s:i(s) =i_v } \mu _s
\rho (s) -\rho (e_{i_v } )) ]=0.
\end{equation}

We define $I_1 =\{ 1\leq v\leq N | s(i_v) =i_u , for\ some\ s\in
\mathcal {R}\}$, and $I_2 = \{ 1\leq v\leq N | s(i_v) \neq i_u ,
for\ any\ s\in \mathcal {R}\}$. There is $\{1,2,\cdots, N \}=I_1
\coprod I_2$.

\begin{equation}
\begin{split}
[ \rho ( e_{i_u } ) , \sum _{v=1} ^{N} ( \sum _{s:i(s) =i_v } \mu _s
\rho (s) -\rho (e_{i_v } )) ] &= -\sum _{v\in I_1  } (e_{i_u }
e_{i_v } -\sum _{s: s(i_v) =i_u } \mu _s e_{i_u } \rho (s) )\\
&+\sum_{v\in I_1 } (e_{i_v } e_{i_u } -\sum _{s: s(i_u ) =i_v } \mu
_s \rho (s)
e_{i_u } )\\
&+ \sum _{v\in I_2 } (e_{i_u } e_{i_v} -e_{i_v } e_{i_u }) \\
&=0.
\end{split}
\end{equation}

Where we used relation (5),(6),(8) in definition 5.1.

\end{proof}

\begin{rem}
Using notations in section 4, it is direct to check the map
$w\mapsto \iota (w) $, $e_i \mapsto p_i $ define a representation
$B_G (\Upsilon) \rightarrow gl(V_G )$. So from $B_G (\Upsilon )$ we
can obtain the generalized Lawrence-Krammer representation.
\end{rem}

\subsection{Low rank cases}

In this subsection we study the algebra $\mathcal {B} _G (\Upsilon
)$ when $G $ is a dihedral group, namely a reflection group of type
$I_2 (m) (m\geq 2)$. These cases are relatively simple but
fundamental. We prove $\mathcal {B} _G (\Upsilon )$ has a cellular
structure and be semisimple for generic parameters.

Denote the dihedral group of type $I_2 (m)$ as $G_m $.

The arrangement of its reflection hyperplanes can be explained with
the following graph.
\begin{figure}[htbp]

  \centering
  \includegraphics[height=4.5cm]{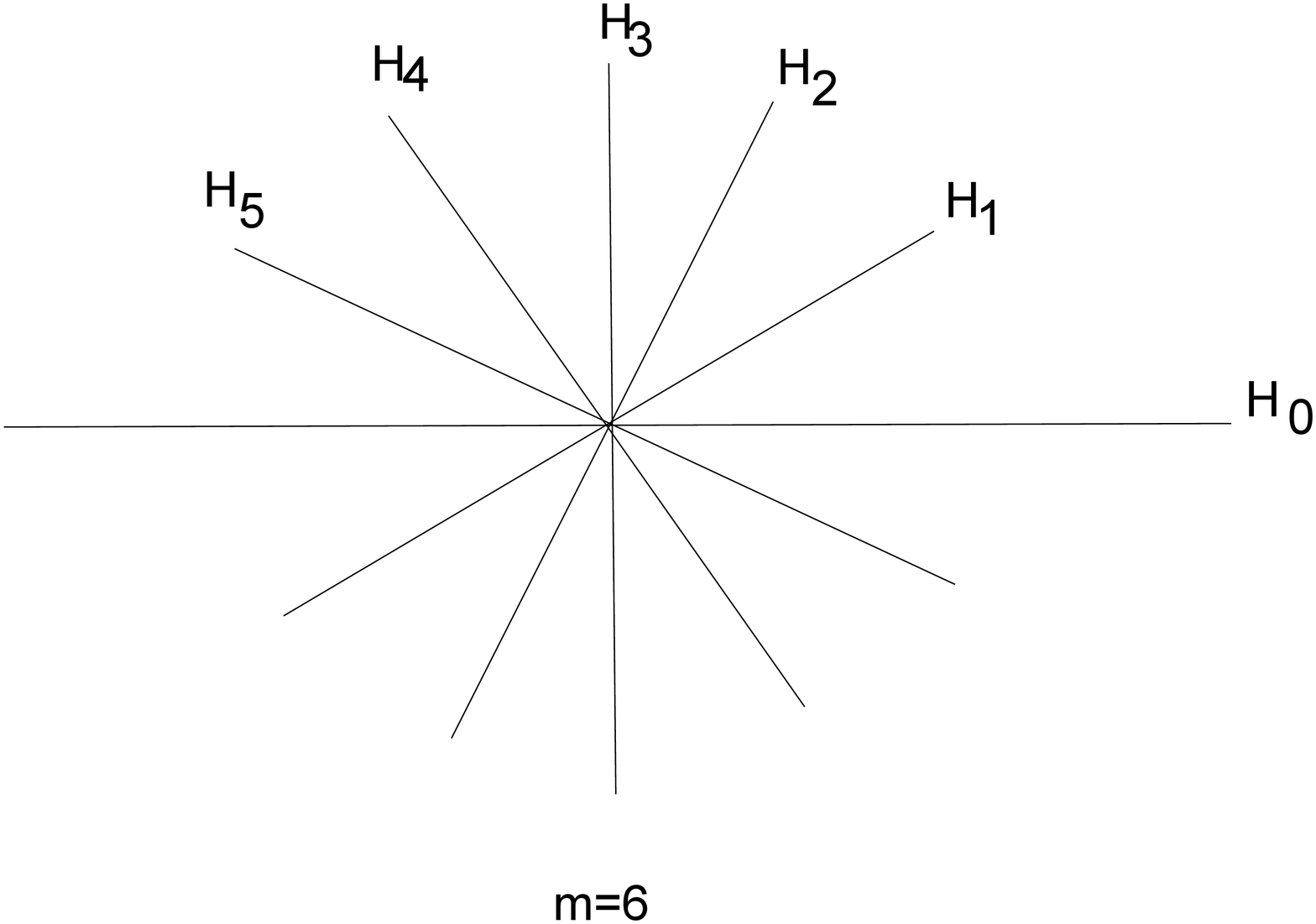}

\end{figure}

There are $m$ lines(hyperplanes) passing the origin. The angle
between every two neighboring line is $\pi /m$. Suppose the $x$-axis
is one of the reflection lines and denote it as $H_0 $, we denote
these lines by $H_0 , H_1 ,\cdots , H_{m-1} $ in anti-clockwise
order as shown in above graph.  Denote the reflection by $H_i $ as
$s_i $. The set of reflections in $G_m$ is $R=\{ s_i \} _{0\leq
i\leq m-1 }$. It is well known that $G_m$ is generated by $s_0 ,s_1
$ with the following presentation
$$< s_0\ ,s_1 , \ | ( s_0 s_1 )^m =1 , s_0 ^2 =s_1 ^2 =1 >.$$

Under this  presentation, $s_i $ can be determined inductively in
the following way.  $s_1 =s_1 $, $s_2 = s_1 s_0 s_1 $ and $s_i =
s_{i-1} s_{i-2} s_{i-1} $. By $[s_i s_j \cdots ]_{k}$ we denote the
unique word starting with $s_i s_j $, in which $s_i $ and $s_j $
appear alternatively, and whose length is $k$. The word $[\cdots s_0
s_1]_k$ is defined similarly. Then $s_i = [s_1 s_0 \cdots ]_{2i-1}
$, and $s_m =[s_1 s_0 \cdots ]_{2m-1 } =s_0 $.

For $k \in \mathbb{Z}$, let $[k]$ be the unique number in
$\{0,1,\cdots , m-1 \}$ such that $k\equiv [k] mod (m \mathbb{Z})$.
It turns out the structure of $B_{G_m}(\Upsilon )$ when $m$ is odd
is rather different from when $m$ is even.\\

\textbf{$B_{G_m}(\Upsilon)$ with $m$ being odd.}  For $0\leq j<i
\leq m-1 $, a number $k(i,j)$ is defined as follows.
  When $i-j$ is even, let $k(i,j) = (i+j)/2$; when $i-j$ is odd, let $k(i,j) = [(m+i+j)/2] $.

Then we have $R(i,j) =\{ k(i,j) \}$ for any $0\leq j<i \leq m-1$.

 For $s_i \in R$, we denote $\mu _{s_i }, \tau _{s_i }$ by $\mu _i
,\tau _i $ simply.  Since all reflections lie in the same conjugacy
class, so all $\mu _i $'s equal some $\mu $ and all $\tau _i$'s
equal some $\tau$. For later use we translate Definition 5.1 in
these cases as follows.

\begin{defi}
 The algebra $B _{G_m } (\Upsilon )$ is generated by $s_0 ,s_1  ; \ e_0 , e_1 ,\cdots ,e_{m-1}$

 with relations:
\begin{enumerate}
\item[(1)] $s_0 ^2 =s_1 ^2 =1 $.
\item[(2)] $(s_0 s_1 )^m =1 $.
\item[(3)] $s_i e_i = e_i s_i = e_i ,\ for\
0\leq i\leq m-1$. Where $s_i = [s_1 s_0 \cdots ]_{2i-1 } $.
\item[(4)] $w e_i  = e_j w $, for\ a word $w$ composed by $s_0 ,\
s_1$ such that \ $ws_i  =s_j w $ is an consequence of relations
$(1)$ and $(2)$.
\item[(5)] $e_i ^2 =\tau e_i ,\ for\ 0\leq i\leq
m-1$.
\item[(6)] $e_i e_j = s_{k(i,j)} e_j = e_i s_{k(i,j)} $ for any $i\neq j$.
\end{enumerate}

\end{defi}

As lemma 5.2 shows, the algebra $B _{G_m } (\Upsilon )$ is spanned
by the set $T_m = G_m \coprod \{ we_i | w\in G_m , 0\leq i\leq m-1
\}$. Because of above relation (3), $| \{ we_i | w\in G_m \}| $
doesn't exceed the number of left cosets of $\{1, s_i \}$ in $G_m $
which equal $m$. So
$$dim B _{G_m } (\Upsilon ) \leq |T_m |\leq 2m+m^2 .$$

We want to show $dim \mathcal {B} _{G_m } (\Upsilon ) \geq 2m+m^2 $.
Let $\pi : B_{G_m } (\Upsilon ) \rightarrow \mathbb{C} G_m $ be the
natural projection. Let $(V_i ,\rho _i )$ $(1\leq i\leq K )$ be all
 irreducible representations of $G_m $. Through $\pi $, every $(V_i
,\rho _i )$ induce a irreducible representation of $B _{G_m }
(\Upsilon )$ which is denoted by $(V_i , \bar {\rho } _i )$. Then
obviously $\bar {\rho } _i $ isn't equivalent to $\bar {\rho } _j$
if $i\neq j$. Beside of these induced representations, $B _{G_m }
(\Upsilon )$ has one more representation coming from Marin's
generalized LK representations. Let $V_{G_m }$ be an vector space
with basis $\{ v_{0 } , \cdots ,v_{m-1 } \}$.  Define a
representation $(V_{G_m}, \rho)$ of $G_m$ by setting $\rho(w)
(v_i)=v_{w(i)}$, where $w(i)$ is determined by $ws_i w^{-1}=
s_{w(i)}$.  For every $i$, define a projection $p_i \in
End(V_{G_m})$ as
$$p_i (v_i)= \tau v_i ; \ p_i (v_j )= \alpha _{i,j } v_i  \ for\ j\neq i,$$
where $\alpha _{i,j} = \# \{k | s_k s_j s_k =s_i \}$. In present
cases $\alpha _{i,j} =1 $ for any $i\neq j$.

\begin{lem}
The map $w\mapsto \rho (w)$, for $w\in G_m $, $e_i \mapsto p_i $
extends to a representation $\rho _{Kr} : B _{G_m } (\Upsilon )
\rightarrow End (V_{G_m })$.

\end{lem}

\begin{proof}
This can be proved by directly checking relations.
\end{proof}

It isn't hard to see that every $v_i $ is a generator of the
representation $(V_{G_m } ,\rho _{Kr })$. Now for $0\neq v= \sum a_i
 v_i \in V_{G_{n}}$, $p_i (v)=( \sum _{j\neq i} a_j + \tau a_i )v_i
.$ We see if $\tau \neq 1 $, there exists some $i$ such that $p_i
(v)\neq 0.$ So we have proved

\begin{lem}
The representation $\rho _{kr}$ is irreducible if $\tau \neq 1$.
\end{lem}

When $\tau \neq 1 $, the algebra $B _{G_m } (\Upsilon )$ have
another irreducible representation $\rho _{Kr } $ whose dimension is
$m$. $\rho _{Kr} $ isn't isomorphic to any $\rho _i $ because $\rho
_{Kr} (e_j ) \neq 0$ but $\rho _i (e_j )=0 $.  By Wedderburn-Artin
theorem,
$$dim B_{G_m } (\Upsilon)\geq \sum _{i=1} ^{K} dim \rho _i ^2 + dim \rho _{Kr}^2
= 2m+m^2.$$

We have proved the following theorem.

\begin{thm}
Suppose $\tau \neq 1$. The set $T_m $ is a basis of $B _{G_m }
(\Upsilon )$. In particular $dim B _{G_m } (\Upsilon )=2m+m^2 $. $B
_{G_m } (\Upsilon )$ is a semisimple algebra.

\end{thm}

It is known there is a cellular structure $(\Lambda , M, C,* )$ on
$\mathbb{C}G_m $. The cellular structure $(\bar {\Lambda },\bar {
M}, \bar {C} ,*)$ of $B _{G_m } (\Upsilon )$ is as follows.

$\bullet$ $\bar {\Lambda } = \Lambda \coprod \{\lambda _{Kr } \}$.We
keep the original partial order in $ \Lambda$, and for any $\lambda
\in \Lambda $, let $\lambda _{Kr} \prec \lambda $.

$\bullet$ For $\lambda \in \Lambda $, $ \bar {M} (\lambda )=
M(\lambda ) $. $\bar {M} (\lambda _{Kr} ) = \{0,1,\cdots ,m-1 \}$.

$\bullet$ For $\lambda \in \Lambda $, and $S,T \in \bar {M}
(\lambda)$ $\bar {C} _{S,T} ^{\lambda }= C_{S,T} ^{\lambda }$. For
$i,j \in \bar {M} (\lambda _{Kr} ) = \{0,1,\cdots ,m-1 \}$, $\bar
{C} _{i,j} ^{\lambda _{Kr }} = we_j $, such that $w(j)=i$. It is
well-defined because $w$ and $ws_j $ are the only elements in $G_m $
satisfying $w^{\prime }(j)=i$. But $ws_j e_j =we_j $.

$\bullet$ Define $*$ be the involution as defined in Lemma 5.5

It isn't hard to prove

\begin{thm}
Above $(\bar {\Lambda },\bar { M}, \bar {C} ,*)$ is a cellular
structure for $B _{G_m } (\Upsilon )$.

\end{thm}

\textbf{$B_{G_m}(\Upsilon)$ with $m$ being even.} The following
facts are true. For $0\leq j<i\leq m-1 $,

$\bullet$ if $i-j$ is even, let $k(i,j) = (i+j) /2$, $k(i,j)
^{\prime} =[(i+j +m)/2]$, then $\{k | s_k s_j s_k = s_i \}=\{k(i,j)
,k(i,j) ^{\prime } \}.$

$\bullet$ if $i-j$ is odd, then $\{k | s_k s_j s_k = s_i
\}=\varnothing $.

The set $R$ of reflections in $G_m $ consists of two conjugacy
classes. $R = R_0 \coprod R_1 $, where $R _0 =\{s_0 ,s_2 ,\cdots
,s_{m-2} \}$ and $R _1 =\{ s_1 ,s_3 ,\cdots , s_{m-1} \}$. In the
datum $\Upsilon$, for convenience we denote $\mu _{s_i } $ as $\mu
_i $, and $\tau _{s_i }$ as $\tau _i $. We have:   $\mu _i =\mu _0 $
, $\tau _i =\tau _0 $ if $i$ is even; $\mu _i =\mu _1$, $\tau _i
=\tau _1 $ if $i$ is odd. So the datum $\Upsilon$ consists of 4
parameters $\mu _0 , \mu _1 , \tau _0 , \tau _1$. An explicit form
of Definition 5.1 in these cases is as follows.

\begin{defi} The algebra $B_{G_m } (\Upsilon )$ is generated by $s_0 ,s_1 ;\  e_0 , e_1 ,\cdots , e_{m-1} $

with the following relations.
\begin{enumerate}
 \item[(1)] $s_0 ^2 =s_1 ^2 =1 $.
 \item[(2)] $(s_0 s_1 )^m =1 $.
\item[(3)] $s_i e_i =e_i s_i $ for all $i$. Where $s_i = [s_1 s_0 \cdots ]_{2i-1}$.
\item[(4)] $w e_i w^{-1 } = e_j $ if $w s_i w^{-1} = s_j $, for $w\in G $.
\item[(5)] $e_i ^2 = \tau _0 e_i$ if $i$ is even, and $e_i ^2 = \tau
_1 e_i $ if $i$ is odd.
\item[(6)] $e_i e_j = 0$ if $i-j $ is odd.
\item[(7)] $e_i e_j = (\mu _{k(i,j)} s_{k(i,j)} + \mu _{k(i,j)^{\prime }}  s_{k(i,j)^{\prime }}) e_j $ if $i-j$ is even.
\end{enumerate}

\end{defi}

Still by lemma 5.2, the algebra $B_{G_m } (\Upsilon )$ is spanned by
the subset $T_m = G_m \coprod \{ we_i | w\in G_m ,\ 0\leq i\leq m-1
\} $. Relation (3) above implies that for every $i$, $|\{ we_i |
w\in G_m  \}|\leq |G_m / \{1,s_i \}| =m$. So we have
$$ dim B_{G_m }(\Upsilon ) \leq |T_m | \leq 2m+ m\cdot m =4l +4l^2 .$$
We use Wedderburn-Artin theorem again to show that for generic
$\Upsilon$, $dim B_{G_m } (\Upsilon) \geq 4l+4l^2 $.\\

Let $\pi :B_{G_m} (\Upsilon)\rightarrow \mathbb{C} G_m $ be the
natural projection sending all $e_i $'s to $0$. Through $\pi $,
every irreducible representation of $G_m $ induce an irreducible
representation of $B_{G_m } (\Upsilon)$. Besides of these induced
representations we have four more different representations
described as follows. Two of them are the infinitesimal version of
the generalized Krammer representation. Let's write down the flat
connection for the generalized Krammer representation in present
case.  The relevant space is $V_{G_m } =\mathbb{C} <v_0 ,\cdots ,
v_{m-1}>$. The action $\rho $ of $G_m $ on $V_{G_m }$ is
$$ \rho (w) (v_i) =v_j
, \ if \ ws_i w^{-1} =s_j \ for \ w\in G_m\ . $$

For every $0\leq i\leq m-1 $, we define a projector $p_i \in End
(V_{G_m })$ as follows.  For $n\in \mathbb{Z}$ define $\epsilon
(n)=0$ if $n$ is even, $\epsilon (n)=1$ if $n$ is odd.

 $p_i (v_i )= \tau _0 v_i $ if $i$ is even; $p_i (v_i ) = \tau _1 v_i
$ if $i$ is odd.

$p_i (v_j )= (\mu _{\epsilon ([\frac {i+j}{2}])} +\mu _{\epsilon
([\frac {i+j}{2} +l ])}) v_i $ if $i-j$ is even.

$p_i (v_j )=0 $ if $i-j$ is odd.

It isn't hard to prove the following lemma.

\begin{lem}
(1)The map $w\mapsto \rho (w) $ for $w\in G_m $; $e_i \mapsto p_i $
for $0\leq i\leq m-1 $ extends to a representation $Kr : B_{G_m }
(\Upsilon) \rightarrow End (V_{G_m })$. \\
(2) Let $V_{G_m } ^0 = \mathbb{C} \langle v_{2i}\rangle _{0\leq
i\leq l-1 }$, and $V_{G_m } ^1 = \mathbb{C} \langle v_{2i+1}\rangle
_{0\leq i\leq l-1 }$. Then $V_{G_m } = V_{G_m } ^0 \oplus V_{G_m }
^1 $ is a decomposition as $B_{G_m } (\Upsilon ) $ representations.
 We denote the subrepresentation on $V_{G_m } ^0 $ as $Kr ^0 $, and
the subrepresentation on $V_{G_m } ^1 $ as $Kr ^1 $.
\end{lem}

Define two matrix $M^0 = ( m^0 _{i,j} )_{l\times l}$ and $M^1 =(m^1
_{i,j} ) _{l\times l} $ whose entries are determined by $p_{2(i-1)}
(v_{2(j-1})) = m^0 _{i,j} (v_{2(i-1)})$ and $p_{2i-1} (v_{2j-1})
=m^1 _{i,j} v_{2i-1}$.

\begin{lem}
If $Det M^0 \neq 0 $ $( Det M^1 \neq 0 )$, then $Kr ^0 $ $(Kr ^1 )$
are different irreducible representations of $B_{G_m } (\Upsilon )
$. Both of them are different from those representations induced
from $\pi $.
\end{lem}

\begin{proof}
For $Kr ^0 $, first every vector $v_{2i }$ is a generator of the
representation since $G_m $ acts on the $\{ v_{2i }\}_{0\leq i\leq
l-1 }$ transitively. Then for any nonzero vector $v= \Sigma _{0\leq
i \leq l-1 } \alpha _i v_{2i}$ in $V_{G_m } ^0 $  , because $Det M^0
\neq 0 $ so there is some $0\leq i\leq l-1 $ such that $p_{2i } (v)
= \lambda v_{2i } \neq 0 $. So all nonzero vectors are generators
and $Kr^0 $ is irreducible. The case of $Kr^1 $ can be proved
similarly.  Because in $Kr^0 $ those elements $p_{2i-1 }$ 's act as
zero but in $Kr^1 $ they don't, so $Kr_0 $ isn't isomorphic to $Kr^1
$. In those representations induced from $\pi $, all the elements
$p_i $'s act as zero. So both $Kr^0 $ and $Kr^1 $ are different from
them.

\end{proof}

There are two more representations $K^0 $ and $K^1 $ as follows.
They arise from the left ideal generated by $(s_l e_0 -e_0)$ and the
left ideal generated by $(s_{l+1} e_1 -e_1 )$ respectively.

We equip an invariant Euclidean metric with $V\cong \mathbb{R} ^2$
on which $G_m $ acts. Then for every $0\leq i\leq 2l-1 $, we chose a
unit vector $v_i $ perpendicular to $H_i $, let $\Phi = \{\pm  v_i
\} _{0\leq i\leq 2l-1 }$. The set $\Phi $ is something like a root
system.

Let $W _{G_m } $ be a complex linear space with a basis $\{ w_{i }\}
_{0\leq i \leq 2l-1 }$. For $0\leq i\leq 2l-1 $, define $\bar {s} _i
\in End(W _{G_m } )$ as follows.

$\bar {s}_i (w_j ) = \epsilon ^i _{j,k } w_k $, if $s_i (v_{[j+l]}
)= \epsilon ^i _{j,k }v_{[k+l ]} .$

\begin{lem}
The correspondence $s_i \mapsto \bar {s} _i $ for $0\leq i\leq 2l-1
$ extends to a representation of $G_m $.
\end{lem}

\begin{pf}
The identity $s_i (v_{[j+l]} )= \epsilon ^i _{j,k }v_{[k+l]} $
implies $\epsilon ^i _{j,k } =\epsilon ^i _{k,j }$, so $\bar {s_i }
^2 =id .$ And since for every pair of $0\leq i,j \leq m-1 $, and
every $0\leq k\leq m-1$, $( s _i s _j )^{m_{i,j}} (v_k ) =v_k $, so
$(\bar {s} _i \bar {s} _j )^{m_{i,j}} (w_k )=w_k $. It follows
$(\bar {s} _i \bar {s} _j )^{m_{i,j}}=id $. So the representation is
well defined.
\end{pf}
\\

For every $0\leq l\leq 2l-1 $, we define $p_i \in End (W _{G_m } )$
as follows.

$p_i (w_j ) =0 $ if $2 \nmid i-j$.

$p_{2i} (v_{2j }) = (\mu _{i+j} \epsilon ^{i+j } _{2j,2i} + \mu
_{[i+j+l] } \epsilon ^{[i+j+l]} _{2j,2i}) v_{2i}$.

$p_{2i-1} (v_{2j-1 }) = (\mu _{i+j-1} \epsilon ^{i+j-1 }
_{2j-1,2i-1} + \mu _{[i+j+l-1] } \epsilon ^{[i+j+l-1]} _{2j-1,2i-1})
v_{2i-1}$.

Define two matrix $A^0 =(a ^0 _{i,j}) _{l\times l}$, $A^1 = (a^1
_{i,j} ) _{l\times l}$ whose entries are determined by $p_{2(i-1)}
(w_{2(j-1)})= a^0 _{i,j} w_{2(i-1)} $, $p_{2i-1} (w_{2j-1}) = a^1
_{i,j} w_{2i-1}$ for $1\leq i,j\leq l$. We have

\begin{thm}
(1) The correspondence $s_i \mapsto \bar {s} _i $, $e_i \mapsto p_i
$ for $0\leq i\leq 2l-1 $ extends to a representation of  $B_{G_m }
(\Upsilon) $ on $W_{G_m }$.\\
(2) Let $W^0 _{G_m }= \mathbb{C} \langle w_{2i} \rangle _{0\leq
i\leq l-1}$, and $W^1 _{G_m }= \mathbb{C} \langle w_{2i+1} \rangle
_{0\leq i\leq l-1}$. Then $W_{G_m } = W^0 _{G_m }\oplus  W^1 _{G_m
}$ is a decomposition of $B_{G_m } (\Upsilon) $ representations. We
denote the representation on $W^0 _{G_m }$ $( W^1 _{G_m })$ as $K^0
 $ $(K^1  )$. \\
(3) If $\det (A^0 )\neq 0 $ and $\det (A^1 )\neq 0 $, then $K^0 $
and $K^1 $ are different irreducible representations of $B_{G_m }
(\Upsilon) $. They are different from $Kr ^0  $ and $Kr^1 $,and be
different from those representations induced by $\pi $ also.

\end{thm}

\begin{pf}
It is straightforward to see (2). For (1) we certify relations in
$B_{G_m } (\Upsilon) $.
\end{pf}
\\

Now by  Wedderburn-Artin theorem, when the data in $\Upsilon$
satisfies the conditions in lemma5.7 and theorem 5.4 (3), we have
$$dim B_{G_m } (\Upsilon) \geq \#\{G_m \} + (dimKr^0  )^2 +(dimKr^1  )^2+(dim K^0 )^2+ (dim K^1 )^2=4l+4l^2 .$$

We have proved the following theorem.

\begin{thm}
When the data in $\Upsilon$ satisfies the conditions in Lemma 5.7
and Theorem 5.4 (3), then $B_{G_{2l} } (\Upsilon)$ is a semisimple
algebra having dimension $4l+4l^2 $. The subset $G_{2l} \coprod \{ w
e_i | w\in G_m , \ 0\leq i\leq 2l-1 \}$ is a basis.
\end{thm}

\textbf{Cellular structure for $B_{G_m } (\Upsilon) $.} Still let
$(\Lambda ,M,C, *)$ be the known cellular structure on $\mathbb{C}
G_{2l}$. We prove that $B_{G_{2l} } (\Upsilon) $ has a cellular
structure $(\bar {\Lambda }, \bar{M}, \bar{C},*)$ as follows.

$\bullet $ $\bar{\Lambda}=\Lambda \cup \{\lambda _{Kr ^0 } , \lambda
_{Kr ^1 }, \lambda _{K^0 }, \lambda _{K^1 } \}$. We keep the partial
order in $\Lambda$, and let $ Kr^i \prec \lambda $, $Kr^1 $ for any
$\lambda \in \Lambda$ and any $i$. Let $K^i \prec Kr^j $ for any
$i,\ j$.

$\bullet $ For $\lambda \in \Lambda $, $\bar{M}(\lambda )
=M(\lambda)$.

$\bullet $ $\bar{M}(\lambda _{Kr ^0 })=\bar{M} (\lambda _{K^0
})=\{0,2,\cdots ,2l-2 \}$.$\bar{M}(\lambda _{Kr ^1
})=\bar{M}(\lambda _{K^1 })=\{1,3,\cdots ,2l-1 \}$.

$\bullet $ Denote the involution defined in Lemma 5.5 extending $s_i
\mapsto s_i $, $e_i \mapsto e_i $ for $0\leq i\leq 2l-1 $ as $*$.

$\bullet $ Denote the including map $\mathbb{C} G_{2l} \rightarrow
B_{G_{2l} } (\Upsilon) $ as $j$. For $\lambda \in \Lambda $, $S,T
\in \bar{M} (\lambda )$, let $\bar{C} ^{\lambda } _{S,T}
=j(C^{\lambda } _{S,T})$.

$\bullet $ Let $[\cdots s_0 s_1 ]_i $ be the word having length $i$,
and in which $s_0 $, $s_1$ appear alternatively. Similarly define
words $[\cdots s_1 s_0 ]_i $.  Let  $\bar {C} ^{\lambda _{K^0 }}
_{2i,0} = [\cdots s_0 s_1 ]_{i} (s_l e_0 -e_0 )$ for $0\leq i \leq
l-1 $.  Let $\bar {C} ^{\lambda _{K^0 }} _{0,2i} =*(\bar {C}
^{\lambda _{K^0 }} _{2i,0})$ for $0\leq i \leq l-1 $. Let $\bar {C}
^{\lambda _{K^0 }} _{2j ,2i} = [\cdots s_0 s_1 ]_{j} (\bar {C}
^{\lambda _{K^0 }} _{0,2i})$ for $0\leq i,j\leq l-1 $.

$\bullet $ Let $\bar{C} ^{\lambda _{K^1 }} _{2i+1 , 1} = [\cdots s_1
s_0 ]_i (s_{l+1} e_1 -e_1 )$ for $0\leq i\leq l-1$. Let $\bar{C}
^{\lambda _{K^1 }} _{1 , 2i+1 } =*(\bar{C} ^{\lambda _{K^1 }} _{2i+1
, 1})$ for any $i$. Let $\bar{C} ^{\lambda _{K^1 }} _{2j+1 , 2i+1 }
=[\cdots s_1 s_0 ]_j (\bar{C} ^{\lambda _{K^1 }} _{1 , 2i+1 })$ for
any $i,j$.

$\bullet $ The group $G_m $ acts on $\{e_{2i} \} _{0\leq i\leq l-1
}$ by conjugacion as $w\cdot(e_{2i}) =we_{2i} w^{-1} $. This action
is transitive. Since $|G_m | =4l $, we know for every pair of $i,j$,
there are two classes in $G_m / \{1,s_{2j} \}$ sending $ e_{2j}$ to
$e_{2i}$. Chose $w_1 ,w_2 $ from each of these two classes, Let
$\bar {C} ^{\lambda _{Kr ^0 }} _{2i,2j } = \frac {1}{2}(w_1 +w_2 )
e_{2j}$. This definition is independent of the choice of $w_1 ,w_2$
since $s_{2j } e_{2j} =e_{2j }$. Similarly consider the conjugating
action of $G_m $ on $\{ e_1 ,\cdots , e_{2l-1 } \}$. There are two
classes in $G_m / \{1,s_{2j-1 } \}$ sending $e_{2j-1 }$ to
$e_{2i-1}$. Chose $w_1, w_2 $ from each of these two classes and let
$\bar {C} ^{\lambda _{Kr ^1 }} _{2i-1,2j-1 } = \frac {1}{2}(w_1 +w_2
) e_{2j-1}$.

The following result can be proved by direct computations.
\begin{prop}
Above data $(\bar {\Lambda }, \bar{M}, \bar{C},*)$ define a cellular
structure on $B_{G_{2l} } (\Upsilon) $.
\end{prop}

\begin{rem}
From results as above it isn't hard to see that the quotient algebra
$\bar{B} _{G_{2l}} (\Upsilon) =B_{G_{2l}} (\Upsilon) / <(s_l e_0
-e_0 ), (s_{l+1} e_1 -e_1 )>$ also satisfy our requests: generically
semisimplity, cellularity, existence of a flat connection, implying
generalized LK representation, but it has simpler structure than
$B_{G_{2l}} (\Upsilon)$. So it can be a candidate for the version
"generalized Brauer algebra " too.  Such quotient algebra as
$\bar{B} _{G_{2l}} (\Upsilon)$ only exist when the Dynkin diagram
contains an even edge, or there is some entry of Coxeter matrix
$m_{i,j} \in 2\mathbb{Z} $. For a general pseudo reflection group
$G$ we can define a quotient algebra $\bar{B} _G (\Upsilon )$ by
adding a relation to Definition 5.1 as:

$(7)$.  $s e_i =s^{'} e_i$  if $s (H_i ) =s^{'} (H_i ) \neq H_i  $
for $s,s^{'} \in R $ and $i\in P$.

\end{rem}

\subsection{Canonical Presentations for Real and G(m,1,n) Cases }

We define algebras $B^{'} _G  (\Upsilon )$ when $G$ is a real or a
cyclotomic reflection group of type $G(m,1,n)$. These algebras have
canonical presentations.  Then we prove $B^{'} _G (\Upsilon )$ is
isomorphic to $B _G (\Upsilon )$. First we do it in cases of
dihedral groups.

\begin{defi}

The algebra $B^{'} _{G } (\Upsilon)$ have the following presentation
when $G$ is the dihedral group of type $I_2 (n )$ where $n=2k+1$ is
an odd number.

$\bullet$ Generators: $S_0 ,S_1 , E_0 ,E_1 .$

$\bullet$ Relations:  $(1)$ $[S_0 S_1 \cdots ]_n =[S_1 S_0 \cdots
]_n $ ; (2)  $S_0 ^2 =S_1 ^2 =1 $;

$(3)$ $S_0 E_0 = E_0 = E_0 S_0 $, $S_1 E_1 = E_1 = E_1 S_1 $ ;

$(4)$ $E_i ^2 = \tau E_i $ for $i=0,1 $;

$(5)$ $E_0 [S_1 S_0 \cdots ]_{2i-1 } E_0 =\mu E_0 $ for $1\leq i\leq
k$;

$(6)$ $E_1 [S_0 S_1 \cdots ]_{2i-1 } E_1 =\mu E_1 $ for $1\leq i\leq
k$;

$(7)$ $[S_0 S_1 \cdots ]_{2k } E_0 = E_1 [S_0 S_1 \cdots ]_{2k } $ ;

\end{defi}

When $G$ is a dihedral group of type $I_2 (2k )$, the element $s_0
s_k =[s_0 s_1 \cdots ]_{2k}$ generate the center. We denote it as
$c$. Denote the set of reflections in $G$ as $R$.

\begin{defi}
The algebra $B^{'} _{G } (\Upsilon)$ have the following presentation
when $G$ is the dihedral group of type $I_2 (n )$ when $n=2k $ is
even. We simply denote $\mu _{s_i }$, $\tau _{s_i }$ as $\mu _i $,
$\tau _i $ respectively.

$\bullet$ Generators: $S_0 ,S_1 ,E_0 ,E_1 $.

$\bullet$ Relations:  $(1)$ $[S_0 S_1 \cdots ]_n =[S_1 S_0 \cdots
]_n ;$  $(2)$ $S_0 ^2 =S_1 ^2 =1 .$

$(3)$ $S_0 E_0 = E_0 S_0 = E_0 $, $S_1 E_1 = E_1 S_1 = E_1 .$

$(4)$ $E_1 w E_0 = E_0 w E_1 =0 $ for any word $w$ composed by $S_0
,\ S_1 .$

$(5)$ $E_{0} [S_1 S_0 \cdots ]_{2i-1 } E_{0} = (\mu _i +\mu _{[i+k]}
c ) E_0 $ for $1\leq i\leq k$ .

$(6)$ $E_1 [S_1 S_0 \cdots ]_{2i-1} E_1 = (\mu _i +\mu _{[i+k]} c )
E_1 $ for $1\leq i\leq k$ .

$(7)$ $[S_1 S_0 \cdots ]_{2k-1 } E_0 = E_0 [S_1 S_0 \cdots ]_{2k-1 }
$;

$(8)$ $[S_0 S_1 \cdots ]_{2k-1 } E_1 =E_1 [S_0 S_1 \cdots ]_{2k-1
}$.

$(9)$ $E_0 ^2 =\tau _{0 } E_0 $, $E_1 ^2 =\tau _{1 } E _1 $.
\end{defi}

\begin{thm}
If $G$ is a dihedral group, then $B_G (\Upsilon )$ is isomorphic to
$B^{'} _G (\Upsilon )$ defined above.

\end{thm}

\begin{proof}

We consider the cases when $G$ is of type $I_2 (2k) $. The cases for
$G$ of type $I_2 (2k+1)$ are similar and easier.

Here we denote the algebra $B_{G } (\Upsilon)$ as $B$,  and denote
the algebra defined in theorem as $B ^{'}$.  Let $j$ be the morphism
from $\mathbb{C} G $ to $B^{'} $ by mapping $s_i \in G $ to $S_i \in
B^{'}$ for $1\leq i\leq n$. Let $\pi $ be the morphism from $B^{'}$
to $\mathbb{C} G $ by mapping $S_i \in B^{'}$ to $s_i \in G $, $E_i
$ to 0. It is easy to check these two morphisms are well defined.
There is $\pi \circ j =id _{\mathbb{C} G }$, which implies $j$ is
injective. For saving notations we denote $j(w)$ as $w$ for $w\in
G$. By Theorem 5.1, there is also a natural injection $l$ from
$\mathbb{C} G$ to $B $. We will also denote $l(w)$ as $w$ for $w\in
G$, and we will identify $S_i $ with $s_i $ when we need.

For $2\leq 2i \leq 2k-2 $, choose any $w \in G $ such that $s_{2i} =
ws_0 w^{-1}$ and let $E_{2i} = w E_0 w^{-1} $. $E_{2i} $ is well
defined with no dependence on choice of $w$.  For example, choose
$w=[S_1 S_0 \cdots ]_{2i-1} $ so $E_{2i} = [S_1 S_0 \cdots ] _{2i-1}
E_0 [S_1 S_0 \cdots ] _{2i-1}$. Similarly for $3\leq 2i-1 \leq 2k-1
$, define $E_{2i-1} = [S_1 S_0 \cdots ]_{2i-2} E_1 [S_1 S_0 \cdots
]_{2i-2} = [S_1 S_0 \cdots ]_{2i-1} E_1 [S_1 S_0 \cdots ]_{2i-1} $.
Define a map $\phi $ from the set of generators of $B$ to $B^{'}$ as
: $\phi (s_0 ) =S_0 $; $\phi (s_1 ) = S_1 $; $\phi (e_i) = E_i $ for
$0\leq i \leq 2k-1 $. Then $\phi $ extends to a morphism from $B$ to
$B^{'}$.  To prove it we only need to certify that $\phi $ keep all
the relations in Definition 5.1. The cases of relation (1), (2) is
straightforward. So $\phi $ can extend to a map from $\mathbb{C} G
\cup \{ e_0 , \cdots , e_{2k-1} \}$.  We prove the case of relation
(7), others are similar. First consider the relation for $e_{2i}
e_0$. We have

$\phi (e_{2i}) \phi ( e_0 ) = E_{2i} E_0 = [s_1 s_0 \cdots ]_{2i-1}
E_0 [s_1 s_0 \cdots ]_{2i-1} E_0 =s_i E_0 s_i E_0 = s_i (\mu _i +
\mu _{[i+k]} c) E_0  \\
= (\mu _i s_i + \mu _{[i+k]} s_{[i+k]})E_{0}
= \phi ( (\mu _{k(2i,0)} s_{k(2i,0)} + \mu _{k(2i,0)^{'} }
s_{k(2i,0) ^{'}} ) e_0  ).$

Since relation (7) is invariant under conjugation of $w\in G$, And
for any $x\in \{ s_0 , s_1 , e_0 , \cdots , e_{2k-1} \} $ any $w\in
G $, we have $\phi (w v w^{-1}) = w \phi (v) w^{-1} $, so the other
relations for $e_{2i} e_{2j} $ reduce to the cases of $e_{2i} e_0 $.

The relations for $e_{2i+1} e_{2j+1} $ are similar.

Then we define a morphism $\psi : B^{'} \rightarrow B$ by extending
the correspondence $S_0 \mapsto s_0 $, $S_1 \mapsto s_1 $, $E_0
\mapsto e_0 $, $E_1 \mapsto e_1 $.  It is easy to check $\psi$ is
well-defined and be the inverse of $\phi$.

\end{proof}

Suppose $G_M $ is a finite Coxeter group with Coxeter matrix $M=
(m_{i,j})_{n\times n} $. As  in  [17],  the condition for $M$ to be
 a coxeter matrix is: $m_{i,i}=1$; $m_{i,j} =m_{j,i} \geq 2 $ for $i\neq j$.
The group $G_M $ has the following presentation:

$\bullet$ Generators: $s_1 ,s_2 ,\cdots ,s_n $;

$\bullet$ Relations:  $[s_i s_j \cdots ]_{m_{i,j}} =[s_j
s_i\cdots]_{m_{i,j}}$ for $i\neq j$;  $s_i ^2 =1$ for any $i$.

$G_M $ can be realized as a group generated by reflections in some
$n$ dimensional linear space naturally,  through the geometric
representation. In the following we identify $G_M $ with it's image
through the geometric representation, and we denote $G_M $ as $G$.

For $w\in G$, any expression $w= s_{i_1} s_{i_2 }\cdots s_{i_r }$
with minimal length is called a reduced form of $w$, and define the
length of $w$ as $l(w) =r$.

For the purpose of proving a key lemma we need to introduce the
Artin group $A_M $ with coxeter matrix $M$. It has the following
presentation.

$\bullet$ Generators: $\sigma _1 , \sigma _2 ,\cdots , \sigma _n $;

$\bullet$ Relations : $[\sigma _i \sigma _j \cdots  ]_{m_{i,j}}
=[\sigma _j \sigma _i \cdots ]_{m_{i,j}}$ for $i\neq j$.

Here we denote $A_M $ as $A$. Let $A^{+} $ be the monoid generated
with the same set of generators and relations. Let $J : A^{+}
\rightarrow A $ be the natural morphism of of monoids.  We have

\begin{thm}

$J$ is an injective map.
\end{thm}

The following theorem is well known.

\begin{thm}
For any $w\in G $, suppose $l(w) =r $ and let $s_{i_1 } \cdots
s_{i_r } $ and $s_{j_1 } \cdots s_{j_r }$ be two reduced forms of
$w$, then in $A^{+}$ we have $\sigma _{i_1} \cdots \sigma _{i_r } =
\sigma _{j_1} \cdots  \sigma _{j_r } $.

\end{thm}

So there is a well defined injective map $\tau : G \rightarrow
A^{+}$ as follows. For $w\in G$, let $s_{i_1 } \cdots s_{i_k } $ be
a reduced form of $w$ and let $\tau (w) = \sigma _{i_1 } \cdots
\sigma _{i_k }$. Denote the natural map from $A^{+}$ to $G$
extending $\sigma _i \mapsto s_i $ as $\pi $.

In $A^{+}$ we denote $b \prec c $ if there is $a \in A^{+}$ such
that $ab=c $. This define a partial order for $A^{+}$. Here is an
important result in Artin group theory.

\begin{thm}{$[7],[14]$}
For $a\in A^{+}$, if $\sigma _i \prec a$, $\sigma _j \prec a$, then
$[\cdots \sigma _j \sigma _i ]_{m_{i,j}} \prec a$.

\end{thm}

Now we can prove the following lemma.

\begin{lem}
Suppose $G$ acts on a set $S$. Suppose there is a subset $\{v_1
,\cdots , v_n \}$ of $S$ such that

(1) If $m_{i,j}=2k+1 $, then $[s_i s_j \cdots ]_{2k} (v_i )=v_j $;

(2) If $m_{i,j} =2k $, then $[s_i s_j \cdots ]_{2k-1} (v_j ) =v_j $;

(3) If $m_{i,j} =2 $, then $s_i (v_j ) =v_j $.

(4) $s_i (v_i ) = v_i  $.

Then an identity $w s_i w^{-1 } = s_j $ in $G$ implies $w (v_i ) =
v_j $.

\end{lem}

\begin{proof}
We prove it by induction on $l(w)$. When $l(w) =0 $ it is evident.
Suppose the lemma is true when $l(w) < k $ and suppose we have an
identity $w s_i w^{-1} =s_j $ where $l(w) =k$.  If $l(w s_i ) = l(w)
-1$, let $w^{'} = ws_i $. Since $w^{'} s_i (w^{'} )^{-1}  =w s_i
w^{-1} = s_j $, by induction we have $w^{'} (v_i ) = v_j $. Which
implies $w(v_i ) =  w^{'} (v_i ) = v_j $ by (4).

Now suppose $l(ws_i ) = l(w) +1 $. Let $s_{i_1 } \cdots s_{i_k }$ be
a reduced form of $w$. We have $s_{i_1 } \cdots s_{i_k } s_i = s_j
s_{i_1 } \cdots s_{i_k }$. Because both sides of the identity are
reduced forms, by theorem5.8 we have $\sigma _{i_1 } \cdots \sigma
_{i_k } \sigma _i = \sigma _j \sigma _{i_1 } \cdots \sigma _{i_k }
=\tau (ws_i)$.

From the condition $l(ws_i ) = l(w) +1$ we know $i_k \neq i$, so by
theorem5.9 we have $[\cdots \sigma _{i_k } \sigma _i ] _{m_{i_k ,i}}
\prec \tau (ws_i )$. So $\tau (ws_i ) = a [\cdots \sigma _{i_k }
\sigma _i ] _{m_{i_k ,i}} $ for some $a \in A^{+}$. Denote $\pi (a)
$ as $w^{'}$ ,and $\pi ([\cdots  \sigma _i \sigma _{i_k } ] _{m_{i_k
,i} -1})$ as $u$. So $w= w^{'} u $. An argument of length shows
$l(w^{'}) = l(w) -l(u) $. There is $u s_i u^{-1} = s_{i_k }$, and by
(1), (2) , (3) we have $u(v_i ) = v_{i_k }$.  So $w^{'} s_{i_k }
(w^{'})^{-1} = w s_i w^{-1} = s_j $. By induction we have $w^{'}
(v_{i_k } ) = v_j $ which implies $w (v_i ) = w^{'} u (v_i ) = v_j
$.

\end{proof}

   Let $R$ be the set of reflections contained in $G$. Let $\Upsilon =
\{ \mu _s ,\tau _s \} _{s\in R } $ be a set of numbers satisfying $
\mu _s = \mu _{s^{'}}$, $\tau _s =\tau _{s^{'}} $ if $s$ and $s^{'}$
lie in the same conjugacy class. We have the following theorem.

\begin{thm}
The algebra $B_{G } (\Upsilon)$ has the following presentation.

$\bullet$ Generators: $s_1 ,\cdots ,s_n ,e_1 ,\cdots ,e_n $.

$\bullet$ Relations: (1)$s_i ^2 =1 $ for any $i$;

(2) $[s_i s_j \cdots ]_{m_{i,j }} =[s_j s_i \cdots ]_{m_{i,j}}$ for
all $i,j$;

(3)  $s_i e_i =e_i =e_i s_i $ for any $i$;

(4)  $e_i ^2 = \tau _i e_i $,  where $\tau _i = \tau _{s_i }$;

(5)   $s_i e_j =e_j s_i $ and  $e_i e_j =e_j e_i $ if $m_{i,j } = 2
$;

(6) $e_i [s_j s_i \cdots ]_{2l-1 } e_i =\mu _{s_{\epsilon }} e_i $
for $1\leq l\leq k$, where $m_{i,j} = 2k+1 $.

 Where $\epsilon= i(j)$if $l$ is odd (even).

(7) $[s_i s_j \cdots ]_{2k } e_i = e_j [s_i s_j \cdots ]_{2k } $ If
$m_{i,j} = 2k+1 $.

(8)  $e_i w e_j =0 $  for any word $w$ composed from $\{s_i ,s_j \}$
If $m_{i,j}= 2k $.

(9)  $e_i [s_j s_i \cdots ]_{2l-1 } e_i = (\mu _s +\mu _{s^{'}}
c_{i,j} ) e_i $ for $1\leq l\leq k$ ,If $m_{i,j}= 2k $. Where
$c_{i,j}= [s_i s_j \cdots ]_{2k}$, $s= [s_j s_i \cdots ]_{2l-1 } $
and $s^{'} = c_{i,j} ^{-1} s $.

(10) $[s_i s_j \cdots ]_{2k-1 } e_j = e_j [s_i s_j \cdots ]_{2k-1 }
$ If $m_{i,j}= 2k $.

\end{thm}

Denote the algebra presented in theorem as $B_{G }^{\prime }
(\Upsilon)$ . The strategy is to construct a morphism from $B_{G
}^{\prime } (\Upsilon)$ to $B_{G } (\Upsilon)$ and a morphism back,
then prove they are the inverse of each other.

In definition 5.1, we first suppose $G $ is the Coxeter group with
Coxeter matrix $M$, then realize it as a reflection group through
the geometric representation $\rho : G\rightarrow GL(V)$. We
identify $G$ with its image in $GL(V)$, denote $\rho (s_i )$ as
$s_i$.  Since $G$ is real, the index set of reflection hyperplanes
$P$ are in one to one correspondence with the set of reflections
$R$. So it is convenient to denote the reflection hyperplane of
$s\in R $ as $H_s $ and write $e_i $ in the definition as $e_s $
when $H_i =H_s $.  Let $\mathcal {A} =\{H_s \}_{s\in R }$ be the set
of reflection hyperplanes. Define $\phi (s_i ) =s_i \in GL(V) $,
$\phi (e_i ) = e_{s_i }$ for $1\leq i\leq n$.

\begin{lem}
$\phi$ extends to a morphism from $B_{G }^{\prime } (\Upsilon)$ to
$B_{G } (\Upsilon)$.
\end{lem}

\begin{pf}
We need to certify that $\phi $ satisfies all relations of $B_{G
}^{\prime } (\Upsilon)$. (1) to (4) are straightforward . When
$m_{i,j} =2 $, The parabolic subgroup $G_{\{i,j\}}$ of $G$ generated
by $s_i ,\ s_j$ is $\mathbb{Z} _2 \times \mathbb{Z} _2 $.  By [15],
the set $\{ H_s \ | H_s \supseteq H_{s_i } \cap H_{s_j } \}$ are in
one to one correspondence with reflections in $G_{\{i,j\}}$. Thus
there are no other reflection hyperplanes containing $ H_{s_i } \cap
H_{s_j }$ except $H_{S_i } $ and $H_{s_j }$. So by definition 5.1
(4), $\phi (e_i ) \phi (e_j ) = e_{s_i } e_{s_j } = e_{s_j } e_{s_i
} = \phi (e_j ) \phi (e_i )$, so relation (5) is satisfied.

When $m_{i,j} =2k+1 $, to certify that $\phi $ satisfies (6) it is
enough to show $e_{s_i } [s_i s_j \cdots ]_{2l-1} e_{s_i } =\mu
_{\epsilon} e_{s_i }$. Denote the reflection $[s_i s_j \cdots
]_{2l-1} $ as $s$. The identity is equivalent to $se_{s_i } s e_{s_i
}= \mu _{\epsilon} s e_{s_i }$. By (3) of definition 5.1, it is
equivalent to $e_{ss_i s} e_{s_i }= \mu _{\epsilon} s e_{s_i }$. By
knowledge of dihedral groups , the only reflection conjugating $s_i
$ to $ss_i s$ is $s$. So the last identity is a consequence of
definition5.1 (5).  Because $[s_i s_j \cdots ]_{2k} s_i [s_i s_j
\cdots ]_{2k}^{-1} = s_j $, so by Definition5.1(3), we have $[s_i
s_j \cdots ]_{2k} e_{s_i } = e_{s_j } [s_i s_j \cdots ]_{2k}$. So
$(7)$ is also satisfied by $\phi$.

When $m_{i,j}=2k >2 $, $s_i $ and $s_j $ don't lie in the same
conjugacy class. let $w$ be any element in the parabolic subgroup of
$G$ generated by $\{s_i ,s_j \}$. Denote $ws_i w^{-1}$ as $s$. The
reflection $s$ lie in different conjugacy classes with $s_j$. Now
$|\{H_{t} \ | H_t \supset H_s \cap H_{s_j \ } \}|= 2k >2 $, so
conditions for definition5.1(9) are fulfilled and we have $e_s
e_{s_j } =0 $ , which is equivalent to $e_{s_i } w e_{s_i }=0 $. So
$\phi $ satisfies (8).

The facts that $\phi $ satisfies (9),(10) can be proved similarly
like the case of  (6), (7).

\end{pf}

Now we construct a morphism from $B_{G } (\Upsilon)$ to $B_{G
}^{\prime } (\Upsilon)$. First we need the following lemma. Let $G$
 acts on $R$ by conjugation.

Through making quotient of $B_{G }^{\prime } (\Upsilon)$ over the
ideal generated by $e_i 's$, it isn't hard to prove that  the
morphism $J$ from $\mathbb{C} G $ to $B_{G }^{\prime } (\Upsilon)$
by sending $s_i $ to $s_i $ is injective.  So for $w\in G$ we
identify $J(w)$ with $w $ .  And denote the imbedding image of $R$
in $B_{G }^{\prime } (\Upsilon)$ as $R^{\prime }$. Set $E^{\prime }
= \{we_i w^{-1} \}_{w\in G, 1\leq i\leq n }$, $E =\{e_s \}_{s\in
R}$. By definition of $B_{G }^{\prime } (\Upsilon)$, the conjugating
action of $G$ on $E^{\prime }$ satisfies conditions in lemma5.10. So
the map $e_i \mapsto e_{s_i }$ $(1\leq i\leq n)$ extends to a
$G$-equivariant bijection $\varphi :E^{\prime } \rightarrow E $.
Define a map $\psi : E\cap R \rightarrow E^{\prime }\cap R^{\prime }
$ by $\psi (e_s ) = \varphi  ^{-1} (e_s )$, $\psi (s)=s$. We have
the following lemma.

\begin{lem}
The map $\psi $ extends to a morphism from $B_{G } (\Upsilon)$ to
$B_{G }^{\prime } (\Upsilon)$. Still denote it as $\psi$.
\end{lem}

\begin{pf}
We prove that $\psi $ satisfies relation (5) in definition5.1. The
fact that $\psi$ satisfies other relations can be proved similarly.
Let $s_{\alpha }, s_{\beta } \in R $ such that $\{ s\in R \ | H_s
\supseteq H_{s_{\alpha }}  \cap H_{s_{\beta }} \} \supsetneq
\{s_{\alpha } ,s_{\beta } \}. $ Suppose the edge $L= H_{s_{\alpha }}
\cap H_{s_{\beta }}$ is in the closure of some Weyl chamber $\Delta
$. Let $\Delta _0 $ be the fundamental Weyl chamber whose walls
consist of $\{H_{s_1 } ,\cdots ,H_{s_n }\}$. Since $G $ acts on the
set of Weyl chambers transitively, there exists $w\in G $ such that
$w \Delta = \Delta _0 $. So, there exists $w\in G$ and $1\leq
i<j\leq n $ such that $w L = H_{s_i } \cap H_{s_j }$.

Denote the parabolic subgroup of $G$ generated by $\{s_i , s_j \}$
as $G_{i,j}$. Since $ws_{\alpha} w^{-1}$ and $ws_{\beta}w^{-1 }$ fix
every point in $H_{s_i } \cap H_{s_j }$,  there exist $\dot{s},
\ddot{s} \in R\cap G_{i,j}$ such that $w\psi (e_{s_{\alpha}} )
w^{-1} =\dot{s}$, $w\psi (e_{s_{\beta}} ) w^{-1} =\ddot{s}$.

Denote the dihedral group of type $I_2 (m_{i,j})$ as $G^{' }$.
Denote the set of reflections in $G^{' }$ as $R^{' }$. Let $\phi $
be the injective morphism from $G^{' }$ to $G$ by sending $s_0$,
$s_1$ to $s_i $,$s_j $ respectively.  We consider the generalized
algebra $B_{G^{'}} (\Upsilon ^{'} )$, where the data $\Upsilon ^{'}
=\{\mu _s , \tau _s \}_{s\in R^{'}}$ is determined by setting $\mu
_s = \mu _{\phi (s)}$, $\tau _s = \tau _{\phi (s)}$ for $s\in R^{'}
$.

By the canonical presentation of $B_{G^{'} } (\Upsilon ^{'} )$ in
theorem5.8 and theorem5.9, the morphism $\phi$ can be extended to a
morphism from $B_{G^{'} } (\Upsilon ^{'} )$ to $B_{G } (\Upsilon)$
by mapping $e_0$, $e_1$ to $e_i $, $e_j $ respectively. We still
denote this morphism as $\phi$. Let $\dot{r}=\phi ^{-1} (\dot{s})$,
$\ddot{r} =\phi ^{-1} (\ddot{s})$. In $B_{G^{'} } (\Upsilon ^{'} )$
we have $e_{\dot{r}} e_{\ddot{r}}=(\Sigma _{r\in R^{'}:
r\dot{r}r=\ddot{r} } \mu _r r ) e_{\ddot{r}} $. Apply $\phi $ to
this identity we get

$$ e_{\dot{s}} e_{\ddot{s}} =( \Sigma _{r\in R^{'}: r\dot{r}r=\ddot{r}
} \mu _r \phi (r) ) e_{\ddot{s}}=( \Sigma _{s\in R
:s\dot{s}s=\ddot{s}} \mu _s s ) e_{\ddot{s}}.$$

The second equality is because $\mu _{\phi (r)} =\mu _r $, and $$\{
\phi (r) \}_{r\in R^{'}: r\dot{r}r=\ddot{r} } = \{s \in R\cap
G_{i,j} \ |\
 s\dot{s}s=\ddot{s} \} =\{s\in R \ |\  s\dot{s}s=\ddot{s} \}.$$

So by conjugating above identity by $w^{-1} $ we have
$$\psi
(e_{s_{\alpha }} ) \psi (e_{s_{\beta}} ) = w^{-1} e_{\dot{s}}
e_{\ddot{s}} w =  ( \Sigma _{s\in R :s\dot{s}s=\ddot{s}} \mu _s
w^{-1}sw ) e_{s_{\beta }}=( \Sigma _{s\in R :ss_{\alpha}s=s_{\beta
}} \mu _s s ) e_{s_{\beta }}.$$

\end{pf}

By definition $\psi$ and $\phi$ are apparently the inverse of each
other, so theorem5.8 is proved. When $M$ is a finite type Coxeter
matrix of simply laced type, any two reflections in $G_M $ are
conjugate to each other. So in the data $\Upsilon$, $\mu _s =\mu $,
$\tau _s = \tau $ for all $s\in R $. Using lemma5.3 we can let $\mu
=1$. So $B_{G_M } (\Upsilon )$ can be also denoted as $B_{G_M }
(\tau )$. Those relations in theorem5.8 in these cases are actually
the same as the relations of generalized Brauer algebras in [10].

\begin{rem}
The algebra $B'_G (\Upsilon )$ can be defined for Coxeter groups of
infinite type as well.
\end{rem}

\textbf{ The Cyclotomic $G(m,1,n)$ cases.} Let $G$ be the cyclotomic
pseudo reflection group of type $G(m,1,n)$. As in [5], let $V $ be a
$n$-dimensional complex linear space with a positive definite
Hermitian metric $<,>$ , let $\{v_1 ,\cdots ,v_n \}$ be a
orthonormal base. Then $G$ can be imbedded in $U(V)$. It's image
consists of monomial matrices whose entries are $m$-th roots of
unit. Here we give a concise description of some facts of $G$
without proof.  Suppose $(z_1 , \cdots ,z_n )$ is the coordinate
system corresponding to $\{v_1 ,\cdots ,v_n \}$. Let $\xi = \exp (
\frac {2\pi \sqrt{-1}}{m})$. For $i\neq j $, $0\leq a\leq m-1 $,
define $H_{i,j;a} = \ker (z_i - \xi ^a z_j )= (v_i -\xi ^a v_j
)^{\perp }$. Define $H_i =\ker (z_i ) =(v_i ) ^{\perp }$. Then $H^a
_{i,j} =H_{j,i;-a}$.

Let $s _{i,j;a }\in U(V)$ be the unique reflection fixing every
points in $H ^a _{i,j} $. Let $s_i $ be the pseudo reflection
defined by: $s_i (v_i )=\xi v_i $; $s_i (v_j ) =v_j$ for $j\neq i$.
Then the set $\mathcal {A}$ of reflection hyperplanes of $G$ is $\{
H _{i,j;a } \} _{i<j ; 0\leq a\leq m-1 } \cup \{ H_i \} _{1\leq
i\leq n }$. The set of pseudo reflections $R$ of $G$ is
$$\{ s _{i,j ;a } \} _{i<j ;0\leq a\leq m-1 } \amalg (\amalg _{k=1} ^{m-1} \{ s^k _i \}_{1\leq i\leq n } ).$$
The left side of above identity gives a decomposition of $R$ into
conjugacy classes. It is well known that $G$ has a canonical
presentation also as follows.

\begin{prop}

If set $S_0 = s_1 $; $S_i = s_{i,i+1 ; 0}$ for $1\leq i\leq n-1 $,
we have the following canonical presentation of $G$.

$\bullet$ Generators:   $ S_0 , \cdots , S_{n-1} $ .

$\bullet$ Relations:    $S_0 ^m = S_i ^2 =id $ for $1\leq i\leq n-1
$; $S_i S_{i+1} S_i =S_{i+1 } S_i S_{i+1} $ for $1\leq i \leq n-2 $;

$S_i S_j = S_j S_i $ for $\mid i-j \mid \geq 2$; $S_0 S_1 S_0 S_1
=S_1 S_0 S_1 S_0 $.

\end{prop}

Now we have a look at the algebra $B_G (\Upsilon) $. The data
$\Upsilon $ now essentially consists of $\mu , \mu _1 , \cdots , \mu
_{m-1} , \tau _0 ,\tau _1$. Where $\mu _{s_{i,j;a }} =\mu $, $\mu
_{s^k _i } = \mu _k $; $\tau _{H _{i,j;a }} = \tau _1 $, $\tau _{H_i
} =\tau _0 $.

The following is a explicit definition of $B_G (\Upsilon )$
according to definition5.1.

\begin{defi} Here is the explicit definition of $B_G (\Upsilon )$
when $G$ is the pseudo reflection group of type $G(m,1,n )$.

$\bullet $ Generators: $\{ \bar{w}  \}_{w\in G } \cup \{ e_{i,j ;a}
\}_{i<j ; 0\leq a\leq m-1 } \cup \{ e_i \} _{1\leq i\leq n }$.

$ \bullet$ Relations: (0) $\bar{w_1 } \bar{w_2 } =\bar{w_3 }$ if
$w_1 w_2 =w_3 $ in $G$;

(1) $\bar{s} _{i,j;a} e_{i,j;a} =e_{i,j;a} \bar{s}_{i,j;a} =
e_{i,j;a}$; $s_i e_i =e_i s_i = e_i $;

(2) $(e_{i,j;a} ) ^2 =\tau _1 e_{i,j;a} $; $(e_i )^2 =\tau _0 e_i
$;

(3) $\bar{w} e_{i,j;a} =e_{k,l;b } \bar{w}$ if $w (H _{i,j;a} )=
H_{k,l;b }$; $\bar{w} e_i =e_j \bar{w} $ if $w(H_i ) =H_j $;

(4) $e_{i,j;a} e_{k,l;b} =e_{k,l;b} e_{i,j;a} $ if $\{i,j \} \cap \{
k,l \} =\emptyset $;  $e_{i,j;a} e_k =e_k e_{i,j;a} $ if $k\not\in
\{i,j \}$;

(5) $e_{i,j;a} e_{i,k;b} = \mu \bar{s} _{j,k;b-a } e_{i,k;b }$ for
$j\neq k$;

$e_i e_j = (\Sigma _{1\leq a\leq m-1} \mu _a \bar{s} _{i,j;a} ) e_j
$;

(6) $e_{i,j;a} e_i = e_i e_{i,j;a} =0$.

\end{defi}

As in the real case, we define the following algebra $B^{'} _G
(\Upsilon )$ and prove it is isomorphic to $B_G (\Upsilon )$.

\begin{defi}
The algebra $B^{'} _G (\Upsilon )$ is defined with the following
generators and relations.

$\bullet$ Generators:  $S_0 , S_1 ,\cdots ,S_{n-1 } $, $E_0 , E_1
\cdots , E_{n-1}$;

$\bullet$ Relations:  (1) $S_0 ^m = S_i ^2 =id$ for $1\leq i\leq n-1
$; $S_0 S_1 S_0 S_1 =S_1 S_0 S_1 S_0 $;

$S_i S_j =S_j S_i $ for $\mid i-j \mid \geq 2$; $S_i S_{i+1} S_i =
S_{i+1} S_i S_{i+1} $ for $1\leq i\leq n-1 $.

 (2) $E_0 ^2 = \tau _0 E_0 $; $E_i ^2 = \tau _1 E_i $ for $1\leq i\leq
n-1$;

(3) $S_1 (S_0 )^i  S_1 E_0 =E_0 S_1 ( S_0 )^i S_1 $ for any $i$;
$S_i E_0 = E_0 S_i $ for $i\geq 1$;

$(S_0 )^i S_1 (S_0 )^i (E_1 )=E_1 (S_0 )^i S_1 (S_0 )^i $ for any
$i$.

$S_i S_{i+1} E_i = E_{i+1} S_i S_{i+1}$ for $i\geq 1$; $S_i E_j =
E_j S_i $ for $\mid i-j \mid \geq 2$.

(4) $S_i E_i = E_i = E_i S_i $.

(5) $E_1 S_0 ^i E_1 = \mu _i E_1 $ for $1\leq i\leq m-1 $; $E_0 S_1
E_0 = (\Sigma _{i=1} ^{m-1} \mu S_1 S_0 ^i S_1 S_0 ^{-i} ) E_0 $.

(6) $E_i E_j =E_j E_i $ for $\mid i-j \mid \geq 2$; $E_i E_{i+1} =
\mu S_i S_{i+1} S_i E_{i+1}$.

(7) $E_0 W E_1 = E_1 W E_0 =0 $ for any word $W$ of $S_i $'s.
\end{defi}

\begin{thm}
The algebra $B^{'} _G (\Upsilon ) $ is isomorphic to $B_G (\Upsilon
) $.

\end{thm}

\begin{proof}

The strategy of proof of this theorem is the same as for the real
case. We construct a morphism $\Phi $ from $B^{'} _G (\Upsilon ) $
to $B_G (\Upsilon ) $ and a morphism $\Psi $ in reverse direction.
Once these morphisms are constructed, it is straight forward to see
they are inverse of each other so the theorem is proved.  The
morphism $\Phi $ is constructed by setting
$$ \Phi (S_0 ) = s_1 ; \Phi (S_i ) = s_{i ,i+1 ;0 }\  for\  i\geq 1   ; \Phi (E_0 ) = e_1 ; \Phi (E_i )= e_{i,i+1 ,0}
 \ for\ i\geq 1.$$

It isn't hard to certify that $\Phi $ satisfying all relations in
definition5.5, so $\Phi $ can extend to a morphism. To define $\Psi
$, the main step is still the definition of $\Psi (e_i ) $ and $\Psi
(e_{i,j;a })$. The following lemma shows there are well defined
elements $F_i $'s and $F_{i,j;a }$'s such that if we set
$$\Psi (w ) = w  , \Psi (e_i ) = F_i , \Psi (e_{i,j;a}) = F_{i,j;a } ,$$
Then $\Psi $ can extend to a morphism from $B_G (\Upsilon )$ to
$B^{'} _G (\Upsilon )$ by certifying that it all relations in
definition5.4. The proof is almost the same as in proof of lemma5.12
so we skep it.

\end{proof}

\textbf{Construction of $F_i $ and $F_{i,j;a }$} \ \ \ \ The
following lemma is similar to lemma5.10.

\begin{lem}
Let $G$ be the pseudo reflection group of type $G(m,1,n)$. Suppose
$G$ acts on a set $S$ and suppose there is a subset $\{v_0 ,v_1
,\cdots , v_{n-1} \}$ such that:

(1) $( S_0 )^i  S_1 ( S_0 )^i  (v_1 ) =v_1 $ for any $i$; $S_1 ( S_0
)^i S_1 (v_0 ) =v_0 $ for any $i$.

(2) $S_i S_{i+1} (v_i ) = v_{i+1 } $ for $i\geq 1$;  $S_{i} S_{i-1}
(v_i ) = v_{i-1 }$ for $i\geq 2$.

(3) $S_i (v_j ) = v_j $ if $\mid i-j \mid \geq 2$.

(4) $S_i (v_i ) =v_i $.

Where $S_i $'s are generators of $G$ as in proposition5.3. For
convenience of presentation here we denote $H_1 $ as $\mathbb{H} _0
$, $H_{i,i+1 ;0 }$ as $\mathbb{H} _i $ for $i\geq 2$.  Then the
identity $w ( \mathbb{H} _i ) = \mathbb{H} _j $ implies $w (v_i )
=v_j $, where $w\in G$.
\end{lem}

\begin{proof}
In this case instead of using Artin monoid we prove it by direct
computation. Let $\bar{v} _i = S_{i-1} \cdots S_1 (v_0 )$ for $1\leq
i\leq n $;  $\bar{v}^a _{i,j} = (S_{j-1} \cdots S_1 S_0 S_1 \cdots
S_{j-1} )^a S_{j-1} \cdots S_{i+1} (v_i )$ for $i<j-1 $; $ \bar{v}
^a _{i,i+1} = (S_i \cdots S_0 \cdots S_i ) ^a (v_i )$. The following
identities show that the set $\{ \bar{v} _i \}_{1\leq i\leq n} \cup
\{ \bar{v} ^a _{i,j } \}_{i<j }$ is closed under the action of $G$,
and the map $J: \mathcal {A} \rightarrow \{ \bar{v} _i \}_{1\leq
i\leq n} \cup \{ \bar{v} ^a _{i,j } \}_{i<j }:$
$$ H_{i,j;a}\mapsto \bar{v}^a _{i,j};\ \ \ H_i \mapsto \bar{v} _i$$
is $G$ equivariant. Thus our theorem is proved.

(a) $S_0 (\bar{v} _i )=\bar{v} _i $.

$S_0 (\bar {v} _i ) = S_{i-1} \cdots S_2 S_0 S_1 (v_0 ) = S_{i-1}
\cdots S_2 S_1 \cdot S_1 S_0 S_1 (v_0 ) = S_{i-1} \cdots S_2 S_1
(v_0 ) = \bar{v} _i $.

(b) $S_0 (\bar{v} ^a _{1,i } )=\bar{v} ^{a-1} _{1,i }$.

$S_0 (\bar{v} ^a _{1,i} )= S_0 (S_{i-1} \cdots S_0 \cdots S_{i-1})^a
S_{i-1} \cdots S_2 (v_1 ) = (S_{i-1} \cdots S_0 \cdots S_{i-1})^a
S_{i-1} \cdots S_2 S_0 (v_1 )$

$= (S_{i-1} \cdots S_0 \cdots S_{i-1})^{a-1} S_{i-1} \cdots S_1 S_0
S_1 S_0 (v_1 ) = (S_{i-1} \cdots S_0 \cdots S_{i-1})^{a-1} S_{i-1}
\cdots S_2 (v_1 ) $

$= \bar{v} ^{a-1} _{1,i}.$

(c) $S_0 (\bar{v} ^a _{i,j})= \bar{v} ^a _{i,j} $ if $i\geq 2$.

$S_0 (\bar{v} ^a _{i,j}) = S_0 (S_{j-1} \cdots S_1 S_0 S_1 \cdots
S_{j-1} )^a S_{j-1} \cdots S_{i+1} (v_i )$

$= (S_{j-1} \cdots S_1 S_0 S_1 \cdots S_{j-1} )^a S_{j-1} \cdots
S_{i+1} S_0 (v_i )= \bar{v} ^a _{i,j}.$

(d) $S_i (\bar{v} _i ) = \bar{v} _{i+1}$ for $i\geq 1$.

$S_i (\bar{v} _i )= S_i S_{i-1} \cdots S_1 (v_0 ) = \bar{v} _{i+1}
.$

(e) $S_i (\bar{v} _{i+1})= \bar{v} _{i}$  for $i\geq 1$.

Equivalent to (d).

(f) $S_i (\bar{v} _j ) = \bar{v} _j $ if $i\neq 0 $ and $j\neq i,i+1
.$

$j\neq i,i+1  \Leftrightarrow i>j \ or\ i<j-1 .$ If $i >j $ then
$S_i (\bar{v} _j )= S_i S_{j-1} \cdots S_1 (v_0 )$

$= S_{j-1} \cdots S_1 S_i (v_0 ) = \bar{v} _j $; If $i< j-1 $ then
$S_i (\bar{v} _j )= S_i S_{j-1} \cdots S_1 (v_0 )$

$= S_{j-1} \cdots S_{i} S_{i+1} S_i \cdots S_1 (v_0 ) =  S_{j-1}
\cdots S_{i+1} S_{i} S_{i+1} \cdots S_1 (v_0 )$

$= S_{j-1} \cdots S_1 S_{i+1} (v_0 ) =\bar{v} _j .$

(g) $S_i (\bar{v} ^a _{i,l} ) = \bar{v} ^a _{i+1 ,l} $ if $l\geq i+2
.$

First we have

$S_i (S_{l-1} \cdots S_0 \cdots S_{l-1} )= S_{l-1} \cdots S_i
S_{i+1} S_i \cdots S_0 \cdots S_{l-1}$

$ =  S_{l-1} \cdots S_{i+1} S_i S_{i+1} \cdots S_0 \cdots S_{l-1} =
S_{l-1} \cdots S_0 \cdots S_{i+1} S_{i} S_{i+1} \cdots S_{l-1} $

$= (S_{l-1} \cdots S_0 \cdots S_{l-1} )S_i .$

So

$S_i (\bar{v} ^a _{i,l} )= S_i (S_{l-1} \cdots S_0 \cdots S_{l-1} )
^a S_{l-1} \cdots S_{i+1} (v_i )$

$= (S_{l-1} \cdots S_0 \cdots S_{l-1} ) ^a S_i S_{l-1}  \cdots
S_{i+1} (v_i ) $

$= (S_{l-1} \cdots S_0 \cdots S_{l-1} ) ^a S_{l-1} \cdots S_{i+2}
S_i S_{i+1} (v_i ) $

$=(S_{l-1} \cdots S_0 \cdots S_{l-1} ) ^a S_{l-1} \cdots S_{i+2}
(v_{i+1} ) = \bar{v} ^a _{i+1 ,l} .$

(h) $S_i (\bar{v} ^a _{l,i}) = \bar{v} ^a _{l, i+1}$ for $l<i $.

$S_i (\bar{v} ^a _{l,i})= S_i (S_{i-1} \cdots S_0 \cdots S_{i-1} )
^a S_{i-1} \cdots S_{l+1}(v_l )$

$ = (S_i \cdots S_0 \cdots S_i ) ^a S_i S_{i-1} \cdots S_{l+1} (v_l
) = \bar{v} ^a _{l,i+1}$.

(i) $S_i (\bar{v} ^a _{k,l} ) = \bar{v} ^a _{k,l} $ if $\{k,l \}
\cap \{ i,i+1 \} = \emptyset .$

 $S_i (\bar{v} ^a _{k,l} )= S_i (S_{l-1} \cdots S_0 \cdots
S_{l-1})^a S_{l-1} \cdots S_{k+1} (v_k )$

\end{proof}

\subsection{ A variation of  $B_G (\Upsilon ) $ }

In definition5.1 of the algebra $B_G (\Upsilon )$, if we replace the
relation (6) with the following weaker relation

$(6) ^{'} $.  $e_i e_j = e_j e_i $, if $\{ k\in I | H_k \supset H_i
\cap H_j \} \neq \{ i,j \}$ and $R(i,j) = \emptyset $,

we obtain an algebra $\hat{ B } _G (\Upsilon )$ larger than $B_G
(\Upsilon )$. It isn't hard to certify that $\hat{ B } _G (\Upsilon
)$ is also finite dimensional, the connection $\Omega $ as in
proposition5.1 is still flat and $G$-invariant, and $\hat{ B } _G
(\Upsilon )$ still contain the generalized Krammer representations.
In the following we explain the reason that taking $B_G (\Upsilon )$
as the generalized Brauer algebra but not $\hat{ B } _G (\Upsilon
)$.

Two important features of the Brauer algebra $\mathcal {B}_n (\tau )
$ are that for generic $\tau $ the algebra $\mathcal {B}_n (\tau ) $
is semisimple and that for all $\tau $ the algebra $\mathcal {B}_n
(\tau ) $ have the same dimension. We hope that generalized Brauer
algebras should also hold these features. The following two lemmas
show that these features don't hold for $\hat{ B } _G (\Upsilon )$.

Here we consider the cases when $G$ are dihedral groups. Notations
are from section 5.2.  In dihedral group $G$ of type $I_2 (2k+1)$,
for any two reflection $s_i ,\ s_j$  the set $R(i,j) \neq \emptyset
$, so in the corresponding algebra $B_G (\Upsilon )$ the condition
(6) doesn't appear.   we only consider the cases when $G$ is a
dihedral group of type $I_2 (2k)$. The main result of this section
is the following proposition.

\begin{prop}
(1) Let $G$ be dihedral group of type $I_2 (2n)$. Let $H_i ,H_j $ be
a pair of reflection hyperplanes such that $R(i,j) =\emptyset $, and
$s_i s_j \neq s_j s_i $. Suppose $s_i s_j s_i =s_k $. In the algebra
$\hat{B} _G (\Upsilon )$, if $\Upsilon$ satisfies the condition:
$\mu _i \pm \mu _{[n+i]} - \tau _k \neq 0$,  then for any
irreducible representation $(V,\rho )$ of $\hat{B} _G (\Upsilon )$
we have $\rho (e_i ) \rho (e_j ) = \rho (e_j ) \rho (e_i ) =0$.

(2) Let $G$ be dihedral group of type $I_2 (2n)$. Let $H_i ,H_j $ be
a pair of reflection hyperplanes such that $R(i,j) =\emptyset $, and
$s_i s_j =s_j s_i $. In the algebra $\hat{B} _G (\Upsilon )$, if
$\Upsilon $ satisfies the condition $\tau _0 (\mu _1 \pm \mu _0 )
\neq 0 $ and $\mu _2 \pm \mu _{[n+2]} -\tau _{4} \neq 0 $, then for
any irreducible representation $(V,\rho )$ of $\hat{ B } _G
(\Upsilon )$ we have $\rho (e_i ) \rho ( e_j ) = \rho (e_j ) \rho
(e_i )=0$.

\end{prop}

\begin{proof}
(1) First we have $e_i e_j = s_i e_i e_j = s_i e_j e_i = s_i e_j s_i
e_i = e_k e_i  $.   The second equality  is by using relation $(6)
^{'}$.

So $e_k (e_i e_j ) = ( e_k )^2 e_i = \tau _k (e_k e_i ) = \tau _k
e_i e_j . $  On the other hand,

$e_k (e_i e_j ) =(e_k e_j ) e_i = (\mu _i s_i + \mu _{[n+i]}
s_{[n+i]}) e_j e_i $.

Now $c= s_i s_{[n+i]}$ is central in $\hat{B} _G (\Upsilon )$, so in
the irreducible representation $(V,\rho )$, $\rho (c)$ is a
constant. Since $c^2 = 1$, so this constant is $1$ or $-1$. When
$\rho (c) =1$, $(\mu _i \rho (s_i ) + \mu _{[n+i]}\rho( s_{[n+i]} ))
\rho (e_j  ) \rho (e_i  ) = (\mu _i \rho (s_i ) + \mu _{[n+i]}\rho(
s_i )) \rho (e_i ) \rho (e_j ) = (\mu _i +\mu _{[n+i]}) \rho (e_j )
\rho (e_i )$. So $(\mu _i +\mu _{[n+i]} -\tau _k ) \rho (e_i ) \rho
(e_j ) =0$ and the condition in statement (1) imply $\rho (e_i )
\rho (e_j) =0$.  Similarly when $\rho (c) =-1$, we have $(\mu _i
-\mu _{[n+i]} -\tau _k)\rho (e_i ) \rho (e_j ) =0$ and the condition
in (1) implies $\rho (e_i ) \rho (e_j ) =0$.

(2) By checking action of reflections on the set of reflection
hyperplanes it isn't hard to see the cases in statement (2) only
happen when $n=2k+1 $ and $j=[i+2k+1]$.   First we consider the case
$i=0 , \ j= 2k+1$. Remember the central element $c= s_i s_{[2k+1 +i
]} = s_0 s_{2k+1}$ for any $i$, so
$$s_i (e_0 e_{2k+1} ) = s_{[2k+1+i]} s_0 s_{2k+1} e_0 e_{2k+1} = s_{[2k+1+i]} s_0 s_{2k+1} e_{2k+1} e_0 = s_{[2k+1+i]} e_{0} e_{2k+1}.$$
Now let (V,$\rho $) be an irreducible representation of $\hat{B} _G
(\Upsilon )$ in which $c$ acts as $id _V $. On the one hand $\rho
[(\mu_0 s_{2k+2} +\mu_1 s_1 ) e_2 s_1 (e_0 e_{2k+1})] =\rho [ (\mu_0
s_{2k+2} +\mu_1 s_1 ) s_1 ( e_0 )^2 e_{2k+1}] = \rho [ \tau_0 (\mu_0
+\mu_1 ) e_0 e_{2k+1 }]$.  Where the second identity is by using
relations of $B_G (\Upsilon )$ and by using above identity.

 On the other hand $\rho
[(\mu_0 s_{2k+2} +\mu_1 s_1 ) e_2 s_1 (e_0 e_{2k+1})] = \rho [e_0
e_2 s_1 (e_0 e_{2k+1})] =\rho [ s_1 e_2 e_0 e_0 e_{2k+1}] =\rho [
\tau_0 s_1 e_2 e_0 e_{2k+1}] =\rho [ \tau _0 s_1 e_2 e_{2k+1} e_0 ]
.$  Since $R(2, 2k+1) =\emptyset$, and $s_2 s_{2k+1} \neq s_{2k+1}
s_2$, so by statement (1) and by conditions in (2), since $s_2
s_{2k+1} s_2 =s_4 $, we have $\rho (e_2 e_{2k+1 }) =0$. So $\tau_0
(\mu_0 +\mu_1 ) \rho (e_0 e_{2k+1 } ) =\tau _0 \rho (s_1 e_2
e_{2k+1} e_0 )=0$. The condition in statement (2) then implies $\rho
(e_0 ) \rho (e_{2k+1 }) =0$.

\end{proof}

\begin{rem}

Canonical presentations of $B_{G} (\Upsilon )$ may help to find
possible deformations of these algebras, i.e, the generalized BMW
algebras. We also hope to find $B_{G} (\Upsilon )$ in some other
settings, for example, find geometric constructions of these
algebras.

\end{rem}

%%%%%%%%%%%%%%%%%%% References %%%%%%%%%%%%%%%%%%%%%%%%%%%%%%%

{\small

\hspace{-0.70cm} {\sc Department of Mathematics}

\nd{\sc University of Science and Technology of China}

\nd {\sc Hefei 230026 China}

\nd {\sc E-mail addresses}: {\sc Zhi Chen} ({\tt
zzzchen@ustc.edu.cn}).

\end{document}